\documentclass[12pt,a4paper]{article}
\usepackage{pstricks}
\newcommand{\VolumeHeader}{}
\newcommand{\VolumeSerial}{LNS}

\newcommand{\ActivityName}{ {\normalsize {\it 
School on Intersection Theory and Moduli
}}\\
{\normalsize {\it 

}}}
\newcommand{\ActivityDate}{ {\normalsize {\it
Trieste, 9-27 September 2002 
}}}

\usepackage{amsmath}
\usepackage{amssymb}
\usepackage{latexsym}
\usepackage{amscd}
\usepackage{xypic}

\newcommand{\LectureHeader}{Moduli spaces of hyperk{\"a}hler manifolds}
\newtheorem{theorem}{Theorem}[section]
\newtheorem{corollary}{Corollary}
\newtheorem{proposition}[theorem]{Proposition}
\newtheorem{rem}[theorem]{Remark}
\newenvironment{remark}{\begin{rem}\rm}{\end{rem}\rm}
\newtheorem{definition}[theorem]{Definition}
\newtheorem{lemma}[theorem]{Lemma}
\newtheorem{ex}[theorem]{Example}
\newenvironment{example}{\begin{ex}\rm}{\end{ex}\rm}
\newtheorem{exs}[theorem]{Examples}
\newenvironment{examples}{\begin{exs}\rm}{\end{exs}\rm}

\newcommand{\prf}{{\em Proof}. }
\newcommand{\qed}{\hspace*{\fill}$\Box$}
\newcommand{\rpfeil}[2]{\stackrel{#2}{\longrightarrow}}
\newcommand{\IR}{{\mathbb R}}
\newcommand{\IP}{{\mathbb P}}
\newcommand{\IQ}{{\mathbb Q}}
\newcommand{\IZ}{{\mathbb Z}}
\newcommand{\IC}{{\mathbb C}}

\newcommand{\Gr}{{\rm Gr}}
\newcommand{\id}{{\rm id}}
\newcommand{\OO}{{\rm O}}
\newcommand{\SO}{{\rm SO}}

\newcommand{\Def}{{\rm Def}}
\newcommand{\re}{{\rm Re}}
\newcommand{\im}{{\rm Im}}
\newcommand{\dual}{^{\scriptstyle\rm v}}

\newcommand{\kc}{{\cal C}}
\newcommand{\km}{{\cal M}}
\newcommand{\kx}{{\cal X}}
\newcommand{\ky}{{\cal Y}}
\newcommand{\kk}{{\cal K}}
\newcommand{\kp}{{\cal P}}
\newcommand{\kt}{{\cal T}}
\newcommand{\epimorph}{\longrightarrow\!\!\!\!\!\!\!\rightarrow}

\begin{document}
\pagestyle{myheadings}
\markboth{\LectureHeader}{\VolumeHeader}
\markright{\VolumeHeader}
\begin{titlepage}
  
\title{Moduli spaces of hyperk{\"a}hler manifolds\\  and mirror symmetry}

\author{Daniel Huybrechts\thanks{huybrech@math.jussieu.fr}
\\[1cm]
{\normalsize
{\it Institut de Math{\'e}matiques Jussieu, Universit{\'e} Paris 7, France.}}
\\[10cm]
{\normalsize {\it Lecture given at the: }}
\\
\ActivityName 
\\
\ActivityDate 
\\[1cm]
{\small \VolumeSerial}
}
\date{}
\maketitle
\thispagestyle{empty}
\end{titlepage}

\baselineskip=14pt
\newpage
\thispagestyle{empty}
\begin{abstract}
These lectures treat some of the basic features of moduli spaces of
hyperk{\"a}hler manifolds and in particular of K3 surfaces. The relation
between the classical moduli spaces and the moduli spaces of conformal
field theories is explained from a purely mathematical point of view.
Recent results on hyperk{\"a}hler manifolds are interpreted in this
context. The second goal is to give a detailed account of mirror symmetry of
K3 surfaces. The general principle, due to Aspinwall and Morrison, and
various special cases (e.g. mirror symmetry for lattice polarized or
elliptic K3 surfaces) are discussed.
\end{abstract}

\vspace{6cm}

{\it Keywords:} Moduli Spaces, K3 Surfaces, Hyperk{\"a}hler Manifolds,
Mirror Symmetry.

\newpage
\thispagestyle{empty}
\tableofcontents

\newpage
\setcounter{page}{1}

\section{Introduction}
Let $(M,g)$ be a compact Riemannian manifold of dimension $2N$ with
holonomy group ${\rm SU}(N)$. If $N>2$ then there is a unique complex
structure $I$ on $M$ such that $g$ is a K{\"a}hler metric
with respect to $I$. For given $g$
and $I$ a symplectic structure on $M$ is naturally defined by the
associated K{\"a}hler form $\omega=g(I(~),~)$. This construction locally
around $g$ yields a decomposition of the moduli space of all Calabi-Yau metrics
on $M$ as $\km^{\rm met}(M)\cong\km^{\rm cpl}(M)\times\km^{\rm khl}(M)$, where
$\km^{\rm cpl}(M)$ is the moduli space of complex structures on $M$ and
$\km^{\rm khl}(M)$ is the moduli space of symplectic structures. Mirror
symmetry in a first approximation predicts for any Calabi-Yau manifold
$(M,g)$ the existence of another Calabi-Yau manifold $(M\dual,g\dual)$
together with an isomorphism $\km^{\rm met}(M)\cong\km^{\rm met}(M\dual)$
which interchanges the two factors of the above decomposition, e.g.\
$\km^{\rm cpl}(M)\cong\km^{\rm khl}(M\dual)$.

The picture has to be modified when we consider the second type of
irreducible Ricci-flat manifolds. If the holonomy group of a
$4n$-dimensional manifold $(M,g)$ is ${\rm Sp}(n)$,
i.e.\ $(M,g)$ is a hyperk{\"a}hler manifold, then the moduli
space of metrics near $g$ does not split into the product of complex
and k{\"a}hler moduli as above, e.g.\ for a given hyperk{\"a}hler metric
there is a whole sphere $S^2$ of complex structures compatible with $g$.
Hence mirror symmetry as formulated above for Calabi-Yau manifolds
needs to be reformulated for hyperk{\"a}hler manifolds. It still
defines an isomorphism between the metric moduli spaces, but
the relation between complex and symplectic structures are more
subtle. Nevertheless, mirror symmetry is supposed to be much simpler for
hyperk{\"a}hler manifolds, as usually the mirror manifold $M\dual$
as a real manifold is $M$ itself.

These notes intend to explain the analogue of the product
decomposition of the moduli space of metrics on a
Calabi-Yau manifold in the hyperk{\"a}hler situation and to show how
mirror symmetry for K3 surfaces, i.e.\ hyperk{\"a}hler manifolds
of dimension $4$, is obtained by the action of a discrete group.

\bigskip

After recalling the main definitions and facts concerning the
complex and metric
structure of these manifolds in Section \ref{basics}
we will soon turn to the global aspects of their moduli spaces.
In Sections \ref{ms} and \ref{ps} we introduce these moduli
spaces as well as the corresponding period domains. The geometric moduli
spaces are studied via maps into the period domains. This will be explained
in Section \ref{periodmaps}. Some of the main results about compact
hyperk{\"a}hler manifolds can be translated into global aspects of these maps.

Compared to other texts (e.g.\ \cite{Periodes})
on moduli spaces of K3 surfaces we will try 
to develop the theory as far as possible for compact hyperk{\"a}hler manifolds
of arbitrary dimension. The second main difference is that
 we also treat the less classical moduli spaces of certain CFTs.
This will be done from a purely mathematical point of view by considering
hyperk{\"a}hler mani\-folds endowed with an additional B-field, i.e.\ a
real cohomology class of degree two.
This will lead to new features starting in Section \ref{Discrete}, where we
let act a certain discrete group on the various moduli spaces. This section
follows papers by Aspinwall, Morrison, and others.
Using this action mirror symmetry of K3 surfaces will be explained in
Section \ref{msK3}. The advantage of this slightly technical approach is that
various versions of mirror symmetry for (e.g.\ lattice polarized or
elliptic) K3 surfaces can be explained by
the same group action. Of course, explaining mirror symmetry
in these terms is only possible for K3 surfaces or hyperk{\"a}hler manifolds. 
Mirror symmetry for general Calabi--Yau manifolds will usually change
the topology.

The text contains little or no original material. The main goal was to
explain global phenomena of moduli spaces of K3 surfaces, or more
generally of compact hyperk{\"a}hler mani\-folds, and to give a
concise introduction into the main constructions
used in establishing  mirror symmetry for K3 surfaces.

We encourage the reader to consult the survey \cite{Periodes} and the
original articles \cite{Aspi,AM}.

\section{Basics}\label{basics}
In this section we collect the basic definitions and facts 
concerning irreducible holomorphic symplectic
manifolds and compact hyperk{\"a}hler
manifolds. Most of the material will be presented without proofs and we shall
refer to other sources for more details (e.g.\ \cite{Beauv,GHJ}).

\begin{definition}
An \emph{irreducible holomorphic symplectic manifold} ({\rm IHS}, for short) is a simply
connected compact K{\"a}hler manifold $X$, such that $H^0
(X,\Omega_X^2)$ is generated by an everywhere non-degenerate
holomorphic two-form $\sigma$.
\end{definition}

Since an IHS is in
particular a compact K{\"a}hler manifold, Hodge decomposition
holds. In degree two it yields
$$\begin{array}{rcl}
H^2(X,\IC)& =&H^{2,0}(X)\oplus H^{1,1}(X)\oplus
H^{0,2}(X)\\
&=&\IC\sigma\oplus H^{1,1}(X)\oplus \IC\bar\sigma.\\
\end{array}$$
The existence of an everywhere non-degenerate two-form $\sigma \in
H^0 (X,\Omega_X^2)$ implies that the manifold has even complex
dimension $\dim_\IC(X)=2n$. Moreover, $\sigma$ induces an
alternating homomorphism $\sigma : \kt_X \to \Omega_X$. Since the
two-form is everywhere non-degenerate, this homomorphism is
bijective. Thus, the tangent bundle and the cotangent bundle of an
IHS are isomorphic.
Moreover, the canonical bundle $K_X=\Omega_X^{2n}$ is trivialized
by the $(2n,0)$-form $\sigma^n$. Thus, an IHS
 has trivial canonical bundle and, therefore,
vanishing first Chern class ${\rm c}_1(X)$.

In dimension two IHS are also called K3 surfaces (K3={\bf K}{\"a}hler, {\bf
K}odaira, {\bf K}ummer). More precisely, by definition
a K3 surface is a compact complex surface with trivial
canonical bundle $K_X$ and such that $H^1(X,{\cal O}_X)=0$. It is a deep fact that any such surface is also
K{\"a}hler \cite{Siu}. Moreover, $H^1(X,{\cal O}_X)=0$ does indeed imply
that such a surface is simply-connected.

\bigskip

Here are the basic examples.

\begin{examples}
{\bf i)} Any smooth quartic hypersurface $X\subset\IP^3$ is a K3 surface,
e.g.\ the Fermat quartic defined by $x_0^4+x_1^4+x_2^4+x_3^4=0$.

{\bf ii)} Let $T=\IC^2/\Gamma$ be a compact two-dimensional complex
torus. The involution $x\mapsto -x$ has 16 fixed points and, thus, the
quotient $T/\pm$ is singular in precisely 16 points. 
Blowing-up those yields a {\it Kummer} surface $X\to T/\pm$, which is a K3
surface containing 16 smooth irreducible rational curves.

{\bf iii)} An elliptic K3 surface is a K3 surface $X$ together with a
surjective morphism $\pi:X\to\IP^1$. The general fibre of $\pi$ is a smooth
elliptic curve.
\end{examples}

It is much harder to construct higher dimensional examples of IHS and all
known examples are constructed by means of K3 surfaces
or two-dimensional complex tori. The list of known examples has been
discussed in length in the lectures of Lehn (see also \cite{GHJ}).

\medskip

So far we have discussed IHS purely from the complex geometric point of
view. However, the most important feature of this type of manifolds is the
existence of a very special metric.

\begin{definition}\label{DEFHK}
A compact oriented Riemannian manifold $(M,g)$
of dimension $4n$ is called
\emph{hyperk{\"a}hler} ({\rm HK}, for short)
if the holonomy group of $g$ equals ${\rm
Sp}(n)$. In this case $g$ is called a hyperk{\"a}hler metric.
\end{definition}

\begin{remark}
If $g$ is a hyperk{\"a}hler metric, then there exist three complex
structures $I$, $J$, and $K$ on $M$, such that $g$ is K{\"a}hler
with respect to all three of them and such that $K=I\circ
J=-J\circ I$. Thus, $I$ is orthogonal with respect to $g$ and the
K{\"a}hler form $\omega_I:=g(I(\, \, ),\, \, )$ is closed
(similarly for $J$ and $K$). Often, this is taken as a definition
of a hyperk{\"a}hler metric. Note that our condition is stronger,
as we not only want the holonomy be contained in ${\rm Sp}(n)$,
but be equal to it.
\end{remark}

\begin{proposition}
Let $(M,g)$ be a HK. Then for any
$(a,b,c)\in\IR^3$ with $a^2+b^2+c^2=1$ the complex manifold
$(M,aI+bJ+cK)$ is an IHS.
\end{proposition}

Thus, for any HK $(M,g)$
there exists a two-sphere
$S^2\subset\IR^3$ of complex structures compatible with the
Riemannian metric $g$.

\begin{remark}
Let $(M,g)$ be a HK. The associated K{\"a}hler forms
$\omega_I$, $\omega_J$, $\omega_K$ span a three-dimensional subspace
$H^2_+(M,g)\subset H^2(M,\IR)$. In fact, this space will always be
considered as a three-dimensional space endowed with the natural
orientation. If $X=(M,I)$, then $H^2_+(M,g)=(H^{2,0}(X)\oplus
H^{0,2}(X))_\IR\oplus\IR\omega_I$, where the orientation is given by the
base $({\rm Re}(\sigma),{\rm Im}(\sigma),\omega_I)$. In order to see this,
one verifies that the holomorphic two-form $\sigma$ on $X=(M,I)$
can be given as $\sigma=\omega_J+i\omega_K$ (cf.\ \cite{GHJ}).
\end{remark}

\begin{definition}\index{K{\"a}hler cone}
Let $X$ be an IHS. The
\emph{K{\"a}hler cone} $ \kk_X\subset H^{1,1}(X,\IR)$ is the open
convex cone of all K{\"a}hler classes on $X$, i.e.\ classes that
can be represented by some K{\"a}hler form.
\end{definition}

 The most important single result on IHS
is the following consequence of the
celebrated theorem of Calabi--Yau:

\begin{theorem}
Let $X$ be an IHS. Then
for any $\alpha\in\kk_X$  there exists a unique hyperk{\"a}hler
metric $g$ on $M$, such that $\alpha=[\omega_I]$ for $\omega_I=g(I(~),~)$.
\end{theorem}

Thus, on any IHS $X$ the K{\"a}hler cone $\kk_X$
parametrizes all possible hyperk{\"a}hler metrics $g$ 
compatible with the given complex structure. Below we will explain how the
K{\"a}hler cone $\kk_X$ can be described as a subset of $H^{1,1}(X)$.

\begin{remark}\label{gleich}
Thus, an IHS $X$ together with a K{\"a}hler class $\alpha\in\kk_X$ is the
same thing as a HK $(M,g)$ together with a compatible complex structure
 $I$. As a short hand, we write $(X,\alpha)=(M,g,I)$ in this case.
\end{remark}

\begin{definition}
The BB(Beauville--Bogomolov)-form of an IHS $X$
is the quadratic form on $H^2(X,\IR)$ given by 
$$q_X(\alpha)=\frac{n}{2}\int_X\alpha^2(\sigma\bar\sigma)^{n-1}+(1-n)
(\int_X\alpha\sigma^{n-1}\bar\sigma^n)(\int_X\alpha\sigma^n\bar\sigma^{n-1}),$$
where $\sigma\in H^{2,0}(X)$
is chosen such that $\int_X(\sigma\bar\sigma)^n=1$
\end{definition}

For any K{\"a}hler class $[\omega]$ we obtain a $q_X$-orthogonal
decomposition $H^2(X,\IR)=(H^{2,0}(X)\oplus
H^{0,2}(X))_\IR\oplus\IR\omega\oplus 
H^{1,1}(X)_\omega$. Here, $H^{1,1}(X)_\omega$
 is the space of $\omega$-primitive real
$(1,1)$-classes. Note that we get a different decomposition for every
K{\"a}hler class $[\omega]\in\kk_X$, but that the quadratic form $q_X$
does not depend on the chosen K{\"a}hler class.

The following proposition collects the main facts about the BB-form $q_X$.

\begin{proposition}\label{BFUJ}
{\rm i)} For any K{\"a}hler class $[\omega]\in\kk_X$ on an IHS $X$
the BB-form $q_X$  is positive definite on
$(H^{2,0}(X)\oplus H^{0,2}(X))_\IR\oplus\IR\omega$ and negative definite on
$H^{1,1}(X)_\omega$.\\
{\rm ii)} There exists a positive real scalar $\lambda_1$ such that
$q_X(\alpha)^n=\lambda_1\cdot\int_X\alpha^{2n}$ for all $\alpha\in H^2(X)$.\\
{\rm iii)} There exists a positive real scalar $\lambda_2$ such  that $\lambda_2
\cdot q_X$
is a primitive integral form on $H^2(X,\IZ)$.\\
{\rm iv)} There exists a positive real scalar $\lambda_3$ such that
$q_X(\alpha)=\lambda_3\cdot\int_X\alpha^2\sqrt{{\rm td}(X)}$ for all $\alpha\in H^2(X)$.
\end{proposition}

After eliminating the denominator of $\sqrt{{\rm td}(X)}$
by multiplying with a universal coefficient $c_n$ that only depends on $n$
we obtain an integral quadratic form $c_n\cdot\int\alpha^2\sqrt{{\rm
    td}(X)}$. In general this form need not be primitive,
but this will be of no importance for us. Moreover, since any IHS has
vanishing odd Chern classes, $\sqrt{{\rm td}(X)}=\sqrt{\hat
  A(X)}$. (Everything that matters here is that $\sqrt{{\rm td}(X)}$ is
purely topological in this case.)
Therefore, in these lectures we will use the following
modified version of the BB-form.

\begin{definition}
The BB-form $q_X$ of an $2n$-dimensional IHS $X$ is given by
$$q_X(\alpha)=c_n\cdot\int_X\alpha^2\sqrt{\hat A
(X)}.$$
\end{definition}

With this definition we see that $q_X$  only depends on the
underlying manifold $M$, i.e.\ for two different hyperk{\"a}hler metrics $g$
and $g'$ and two compatible complex structures $I$ resp. $I'$ the BB-forms
with the above definition of $X=(M,I)$ and $X'=(M,I')$ coincide.

\medskip

Note for $n=1$ we have $c_1=1$ and thus
$q_X$ is nothing but the intersection pairing $\alpha\cup\alpha$ of the
four-manifold underlying a K3 surface. The quadratic form in this case is
even, unimodular and indefinite and can thus be explicitly determined:

\begin{proposition}
The intersection form $(H^2(X,\IZ),\cup)$ of a K3 surface $X$ is isomorphic
to the K3 lattice $2(-E_8)\oplus 3U$, where $U$ is the standard hyperbolic
plane $\left( \IZ^2 ,
\bigl(\begin{smallmatrix}0&1\\1&0\end{smallmatrix}\bigr)\right)$.
\end{proposition}

\begin{definition}
The BB-volume of a HK $(M,g)$ is
$$q(M,g):=q_X([\omega_I]),$$
where $X=(M,I)$ is the IHS associated to one of the compatible complex
structures $I$ and $\omega_I$ is the induced K{\"a}hler form.
\end{definition}

Note that the BB-volume does not depend on the chosen complex structure.
Analogously one can define the volume of an IHS endowed with a K{\"a}hler
class $\alpha$ as $q_X(\alpha)$.
For a K3 surface one has $q(M,g)=\int\omega_I^2$,
which is the usual volume
up to the scalar factor $1/2$. In higher dimension the usual volume is
of degree $2n$ and the BB-volume is quadratic.
Of course, due to Proposition \ref{BFUJ} one knows that up to a scalar factor
$q(M,g)^{n}$ equals the standard volume, but this factor might a priori depend on
the topology of $M$. 

\bigskip

What makes the theory of K3 surfaces and higher-dimensional HK so pleasant
is that they can be studied  by means of their period.

\begin{definition}
Let $X$ be an IHS. The \emph{period}
of $X$ is the lattice $(H^2(X,\IZ),q_X)$ 
endowed with the weight-two Hodge structure
$H^2(X,\IZ)\otimes\IC=H^2(X,\IC)=\IC\sigma\oplus H^{1,1}(X,\IC)\oplus\IC\bar\sigma$.
\end{definition}

Since $H^{1,1}(X,\IC)$ is orthogonal with respect to $q_X$ and 
$\IC\bar\sigma$ is the complex conjugate of $\IC\sigma$, the period of the
IHS $X$ is in fact given by the lattice $(H^2(X,\IZ),q_X)$ and the line
$\IC\sigma\subset H^2(X,\IC)$.

The theory of K3 surfaces is crowned by the so called Global Torelli
Theorem (due to Pjateckii-Sapiro, Shafarevich,
Burns, Rapoport, Looijenga, Peters, Friedman):

\begin{theorem}\label{GTK3}
Let $X$ and $X'$ be two K3 surfaces and let $\varphi:H^2(X,\IZ)\cong
H^2(X',\IZ)$ be an isomorphism of their periods such that
$\varphi(\kk_X)\cap\kk_{X'}\ne\emptyset$. Then there exists a unique
isomorphism $f:X'\cong X$ such that $f^*=\varphi$.
\end{theorem}

Moreover, an arbitrary isomorphism between the periods of two K3 surfaces
is in general not induced by an isomorphism of the K3 surfaces, but the
K3 surfaces are nevertheless isomorphic.

The uniqueness assertion in the Global Torelli Theorem is  roughly
proven as follows (cf.\ \cite{LP}):
If $f$ is a non-trivial automorphism of finite order with $f^*={\rm id}$ then the holomorphic two-form $\sigma$ is invariant under $f$ and the action
at the fixed points is locally of the form
$(u,v)\mapsto (\xi\cdot u,\xi^{-1}\cdot v)$. Using Lefschetz fixed point formula
and again $f^*={\rm id}$ one finds that there are 24 fixed points. Thus, the
minimal resolution $\tilde X$ of the quotient $X/\langle f\rangle$ contains
24 pairwise disjoint curves. Moreover, one verifies that $\tilde X$ is again
a K3 surface. The last two statements together yield a contradiction.

The Global Torelli Theorem in the above version fails completely
in higher dimensions.
E.g.\ if $f:X\cong X$ is an automorphism of a K3 surface $X$ such that
$f^*={\rm id}$, then $f={\rm id}$. This does not hold in higher dimensions
\cite{Beauville2}. Even worse, due to a recent counterexample of Namikawa
\cite{Namikawa} one knows that higher dimensional IHS $X$ and $X'$ might have
isomorphic periods without even being birational.
A possible formulation of the Global Torelli Theorem in higher
dimensions using derived catgeories was proposed in
\cite{Namikawa}. However,
uniqueness is not expected. Compare the discussion in Section \ref{derGT}.
\medskip

Often, a certain type of K3 surfaces is distinguished by the form of the
period. We explain this in the three examples presented earlier. In fact,
the proofs of these descriptions are all quite involved.

\bigskip

\begin{example}
{\bf i)} Let $X$ be a K3 surface such that ${\rm Pic}(X)=H^2(X,\IZ)\cap
H^{1,1}(X)$ is generated by a class $\alpha$ with $\alpha^2=4$. Then $X$ is
isomorphic to a quartic hypersurface in $\IP^3$ and $\alpha$ corresponds to
${\cal O}(1)$ (cf.\ \cite[Exp. VI]{Periodes}).

{\bf ii)} Let $X$ be a K3 surface such that ${\rm Pic}(X)$ contains 16 
disjoint smooth irreducible rational curves $C_1,\ldots,C_{16}\subset X$
such that $\sum [C_i]\in H^2(X,\IZ)$ is two-divisible. Then $X$ is isomorphic to a Kummer surface.

 This description of Kummer surfaces is not entirely in terms of the
 period.
Later we will rather use the following description of an even more
special type of K3 surfaces: Let $X$ be a K3 surface such that the lattice
$(H^{2,0}(X)\oplus H^{0,2}(X))_\IZ$ is of rank two and any vector $x$
in this lattice satisfies $x^2\equiv0\mod4$. Then $X$ is a Kummer
surface. It turns out that K3 surfaces with this type of period
are exactly the exceptional Kummer surfaces, i.e.\ Kummer surfaces with
${\rm rk}({\rm Pic}(X))=20$ (cf.\ \cite[Exp.VIII]{Periodes}).

{\bf iii)} Let $X$ be a K3 surface such that there exists a class
$\alpha\in H^2(X,\IZ)\cap H^{1,1}(X)$ with $\alpha^2=0$. Then $X$ is an
elliptic K3 surface. Clearly, if $X\to\IP^1$ is an elliptic K3 surface then
the class of the fibre defines such a class. But note that conversely
not every class $\alpha$ with $\alpha^2=0$ is automatically
a fibre class of some elliptic fibration, but by applying
certain reflections it can be be made into one (cf.\ \cite{Beauvillesurf}).
\end{example}

In order to get a better feeling for the set of all possible hyperk{\"a}hler
structures on an IHS $X$ we shall discuss the K{\"a}hler cone in some more
detail.

\begin{definition}
The \emph{positive cone} $\kc_X$ of an IHS $X$ is the connected component
of the open set
$\{\alpha~|~q_X(\alpha)>0\}\subset H^{1,1}(X,\IR)$ that contains the
K{\"a}hler cone $\kk_X$.
\end{definition}

(Here we use the fact that $q_X(\alpha)>0$ for any K{\"a}hler class
$\alpha$.)
Thus, $\kc_X\cup(-\kc_X)$ can be entirely read off the period of $X$. This
is no longer possible for the K{\"a}hler cone, but one can at least try to
find a minimal set of further geometric information that determines $\kk_X$
as an open subcone of $\kc_X$.

\begin{proposition}
The K{\"a}hler cone $\kk_X\subset\kc_X$ is the open subset of all
$\alpha\in\kc_X$ such that $\int_C\alpha>0$ for all rational curves
$C\subset X$. If $X$ is a K3 surface it suffices to test smooth rational
curves (cf.\ \cite{Periodes,BPV,GHJ}).
\end{proposition}

Since any smooth irreducible rational curve $C$ in a K3 surface $X$ defines
a $(-2)$-class $[C]\in H^{1,1}(X,\IZ)$, one can use this result to show
that
for any class $\alpha\in\kc_X$ there
exists a finite number of smooth rational
curves $C_1,\ldots,C_k\subset X$ such that $s_{C_1}\ldots
s_{C_k}(\alpha)\in\overline{\kk}_X$, where $s_{C}$ is the reflection in the
hyperplane $[C]^\perp$. Of course, these reflections $s_C$ are contained in
the discrete orthogonal group $\OO(\Gamma)$ of the lattice
$\Gamma=(H^2(X,\IZ),\cup)$.

\section{Moduli spaces}\label{ms}

Ultimately, we will be interested in moduli spaces of irreducible
holomorphic symplectic manifolds (IHS), hyperk{\"a}hler manifolds (HK), etc.
In this section we will introduce moduli spaces of such manifolds
endowed with an additional marking. A marking in general refers
to an isomorphism of the second cohomology with  a fixed lattice.
The choice of such an isomorphism gives rise to the action of a
discrete group and the quotients by this group will eventually
yield the true moduli spaces. For this section we fix a lattice
$\Gamma$ of signature $(3,b-3)$ and an integer $n$.

\subsection{Moduli spaces of marked IHS}

\begin{definition}
A \emph{marked IHS} is a pair $(X,\varphi)$ consisting of an IHS of
complex dimension $2n$ and a lattice isomorphism
$\varphi:(H^2(X,\IZ),q_X)\cong\Gamma$. We say that two marked IHS
$(X,\varphi)$ and $(X',\varphi')$ are equivalent,
$(X,\varphi)\sim(X',\varphi')$, if there exists an isomorphism
$f:X\cong X'$ of complex manifolds such that
$\varphi'=\varphi\circ f^*$.
\end{definition}

\begin{definition}
The moduli space of marked IHS is the space
$$\kt_\Gamma^{\rm cpl}:=\{(X,\varphi)={\rm marked~IHS}\}/\!\!\sim.$$
\end{definition}

A priori, $\kt^{\rm cpl}_\Gamma$ is just a set, but, as we will see
later, it can be endowed with the structure of a topological space
locally isomorphic to a complex manifold of dimension $b-2$.

Let $X$ be an IHS and $\varphi$ a marking of $X$. If $\kx\to\Def(X)$
is the universal deformation of $X=\kx_0$, then $\Def(X)$ is a
smooth germ of dimension $h^1(X,\kt_X)$. We may represent
$\Def(X)$ by a small disc in $\IC^{h^1(X,\kt_X)}$. The marking
$\varphi$ induces in a canonical way a marking $\varphi_t$ of the
fibre $\kx_t$ for any $t\in\Def(X)$. Using the Local Torelli Theorem
(cf.\ Section \ref{periodmaps}) we see that the induced map $\Def(X)\to\kt^{\rm
cpl}_\Gamma$ is injective, i.e.\ any two fibres of the family
$\kx\to\Def(X)$ define non-equivalent marked IHS. The various
$\Def(X)\subset\kt^{\rm cpl}_\Gamma$ for all possible choices of
$X$ and markings $\varphi$ cover the moduli space $\kt^{\rm
cpl}_\Gamma$. Since the universal deformation $\kx\to\Def(X)$ of
$X=\kx_0$ is, at the same time, also the universal deformation of
all its fibres $\kx_t$, one can define a natural topology on
$\kt^{\rm cpl}_\Gamma$ by gluing the complex manifolds $\Def(X)$.
Thus, locally $\kt^{\rm cpl}_\Gamma$ is a smooth complex manifold
of dimension $h^1(X,\kt_X)=b-2$. However, $\kt^{\rm cpl}_\Gamma$
is not a complex manifold, as it does not need to be Hausdorff. In
fact, no example is known, where $\kt^{\rm cpl}_\Gamma$
would be Hausdorff and conjecturally this never happens.

A family $(\kx,\varphi)\to S$ of marked IHS is a family $\kx\to
S$ of IHS of dimension $2n$ and a family of markings $\varphi_t$
of the fibres $\kx_t$ locally constant with respect to
$t$.

The universality of $\kx\to\Def(X)$ immediately implies the following

\begin{lemma}
If $(\kx,\varphi)\to S$ is a family of marked IHS, then there
exists a canonical holomorphic map $\eta:S\to\kt^{\rm
cpl}_\Gamma$, such that $\eta(t)=[(\kx_t,\varphi_t)]$.\qed
\end{lemma}

\begin{remark}
In order to construct a universal family over $\kt^{\rm
cpl}_\Gamma$ one would need to glue universal families
$\kx\to\Def(X)$, $\ky\to \Def(Y)$, where $(X,\varphi)$ and
$(Y,\psi)$ are marked IHS, over the intersection $\Def(X)\cap
\Def(Y)\subset\kt^{\rm cpl}_\Gamma$. This is only possible if for
$t\in\Def(X)\cap\Def(Y)$ there exists a unique isomorphism
$f:\kx_t\cong\ky_t$ with $\psi_t=\varphi_t\circ f^*$. For K3
surfaces the uniqueness can be ensured due to the strong version
of the Global Torelli Theorem (see Thm.\ \ref{GTK3}),
but in higher dimensions this
fails. Thus, $\kt^{\rm cpl}_\Gamma$ is, in general, only a
coarse moduli space.
\end{remark}


\subsection{Moduli spaces of marked HK}

\begin{definition}
A \emph{marked HK} is a triple $(M,g,\varphi)$, where $(M,g)$ is a
compact HK of dimension $4n$ in the sense of Proposition \ref{DEFHK} and
$\varphi$ is an isomorphism $(H^2(M,\IZ),q)\cong\Gamma$. Two
triples $(M,g,\varphi)$, $(M',g',\varphi')$ are equivalent,
$(M,g,\varphi)\sim(M',g',\varphi')$, if there exists an isometry
$f:(M,g)\cong(M',g')$ with $\varphi'=\varphi\circ f^*$.
\end{definition}

\begin{definition}
The moduli space of marked HK is the space
$$\kt^{\rm met}_\Gamma:=\{(M,g,\varphi)={\rm marked~HK}\}/\!\!\sim.$$
\end{definition}

A slightly different approach towards $\kt_\Gamma^{\rm met}$
will be explained in Section \ref{without}. There, the manifold
$M$ is fixed and only the metric $g$ is allowed to vary.

\subsection{Moduli spaces of marked complex HK or K{\"a}hler IHS}

Recall (cf.\ Remark \ref{gleich}) that there is a bijection
between  HK with a compatible complex structure and IHS with a
chosen K{\"a}hler class. Thus, the two moduli spaces are naturally
equivalent.

\begin{definition}
A \emph{marked complex HK} is a tuple $(M,g,I,\varphi)$, where
$(M,g,\varphi)$ is a marked HK and $I$ is a compatible complex
structure on $(M,g)$. A \emph{marked K{\"a}hler IHS} is a triple
$(X,\alpha,\varphi)$, where $(X,\varphi)$ is a marked IHS and
$\alpha\in\kk_X$ is a K{\"a}hler class.
Two marked complex HK $(M,g,I,\varphi)$,
$(M',g',I',\varphi')$ are equivalent if  there exists an isometry
$f:(M,g)\cong(M',g')$ with $I=f^*I'$ and
$\varphi'=\varphi\circ f^*$. Analogously, one defines the
equivalence of marked K{\"a}hler IHS.
\end{definition}

Note that the equivalence relation is compatible with the natural
bijection $\{(M,g,I,\varphi)\}\leftrightarrow
\{(X,\alpha,\varphi)\}$.

\begin{definition}
The moduli space of complex HK or, equivalently, of K{\"a}hler IHS
is the space
\begin{eqnarray*}
\kt_\Gamma&:=&\{(M,g,I,\varphi)={\rm marked ~complex ~HK}\}/\!\!\sim\\
&=&\{(X,\alpha,\varphi)={\rm
marked~K\ddot{a}hler~IHS}\}/\!\!\sim.\\
\end{eqnarray*}
\end{definition}

Obviously, there are two forgetful maps $m:(M,g,I,\varphi)\mapsto
(M,g,\varphi)$ and $c:(X,\alpha,\varphi)\mapsto (X,\varphi)$. The
following diagram is the hyperk{\"a}hler version of the product
decomposition of the metric moduli space for Calabi--Yau manifolds.

$$\xymatrix{
\kt_\Gamma\ar[r]^{m}\ar[d]^{c}&\kt^{\rm met}_\Gamma\\
\kt^{\rm cpl}_\Gamma&}$$

\begin{proposition}\label{twistorlines}
The set $\kt_\Gamma$ has the structure of a real manifold of
dimension $3(b-2)$. The fibre $c^{-1}(X,\varphi)=\kk_X$ is a real
manifold of dimension $b-2$. The fibre $m^{-1}(M,g,\varphi)$ is
naturally isomorphic to the complex manifold $\IP^1$. The induced
map $c:\IP^1=m^{-1}(M,g,\varphi)\to\kt^{\rm cpl}_\Gamma$ is a
holomorphic embedding. The map $m:c^{-1}(X,\varphi)\to\kt^{\rm
met}_\Gamma$ is a real embedding.\qed
\end{proposition}

The line $\IP^1\subset\kt_\Gamma^{\rm cpl}$ is also called `twistor line'.
Disposing of a global deformation like this, is one of the key tools in studying
moduli spaces of IHS.

\subsection{CFT moduli spaces of HK}

From a geometric point of view the following moduli space is an
almost trivial extension of $\kt_\Gamma$. However, it will become
of central interest in later sections, when we will let act the
full modular group on it. This group action will relate very
different HK and thus gives rise to mirror symmetry phenomena.

\begin{definition}
A marked complex HK with a B-field is a tuple
$(M,g,I,B,\varphi)$, where $(M,g,I,\varphi)$ is a marked complex
HK and $B\in H^2(M,\IR)$. Two such tuples $(M,g,I,B,\varphi)$,
$(M',g',I',B',\varphi')$ are equivalent if there exists an
isometry  $f:(M,g)\cong(M',g')$ with $I=f^*I'$,
$\varphi'=\varphi\circ f^*$, and $f^*(B')=B$.
\end{definition}

\begin{definition}
The $(2,2)$-CFT moduli space of HK is the space
$$\kt_\Gamma^{\scriptscriptstyle{(2,2)}}:=\{(M,g,I,B,\varphi)={\rm marked ~complex ~HK
~with ~B{\rm -}field}\}/\!\!\sim.$$
\end{definition}

Clearly, the moduli space $\kt_\Gamma^{\rm
\scriptscriptstyle{(2,2)}}$ is naturally isomorphic to
$\kt_\Gamma\times\Gamma\otimes\IR$ by mapping $(M,g,I,B,\varphi)$
to $((M,g,I,\varphi),\varphi_\IR(B))$. In particular,
$\kt_\Gamma^{\rm \scriptscriptstyle{(2,2)}}$ is  a real manifold
of dimension $4b-6$.

\medskip

Analogously, one defines the $(4,4)$-CFT moduli space
$$\kt_\Gamma^{\scriptscriptstyle{(4,4)}}:=\{(M,g,B,\varphi)={\rm marked ~HK
~with ~B{\rm -}field}\}/\!\!\sim.$$
In particular, there is a natural
forgetful map $\kt_\Gamma^{(2,2)}\to\kt_\Gamma^{(4,4)}$
which is surjective with fibre $S^2$.
\subsection{Moduli spaces without markings}\label{without}

All previous moduli spaces parametrize various geometric objects
with an additional marking of the second cohomology. Of course,
what we are really interested in are the true moduli spaces
$\km_\Gamma^{\rm cpl}$, $\km_\Gamma^{\rm met}$, $\km_\Gamma$,
$\km^{(2,2)}_\Gamma$, and $\km^{(4,4)}$. E.g.\ $\km^{\rm
cpl}_\Gamma$ is the moduli space of  IHS $X$ of dimension $2n$
such that $(H^2(X,\IZ),q_X)$ is isomorphic to $\Gamma$, but
without actually fixing the isomorphism. Analogously for the
other spaces. In other words one has:
\begin{eqnarray*}
\km_\Gamma^{\rm
cpl}=\raisebox{-.6ex}{$\OO(\Gamma)$}\setminus\raisebox{.3ex}{$\kt_\Gamma^{\rm
cpl}$} ,& ~~&\km_\Gamma^{\rm
met}=\raisebox{-.6ex}{$\OO(\Gamma)$}\setminus\raisebox{.3ex}{$\kt_\Gamma^{\rm
met}$}, ~~\km_\Gamma^{\rm
}=\raisebox{-.6ex}{$\OO(\Gamma)$}\setminus\raisebox{.3ex}{$\kt_\Gamma^{\rm
}$},\\
\km_\Gamma^{\rm
(2,2)}=\raisebox{-.6ex}{$\OO(\Gamma)$}\setminus\raisebox{.3ex}{$\kt_\Gamma^{\rm
(2,2)}$} ,& ~~&\km_\Gamma^{\rm
(4,4)}=\raisebox{-.6ex}{$\OO(\Gamma)$}\setminus\raisebox{.3ex}{$\kt_\Gamma^{\rm
(4,4)}$}.\\
\end{eqnarray*}

The Teichm{\"u}ller spaces $\kt_\Gamma^*$ are in general better
behaved. E.g.\ the moduli spaces are usually singular at points
that correspond to manifolds with a bigger automorphism group
than expected. This usually leads to orbifold singularities.
However, sometimes the passage from the Teichm{\"u}ller space to
the moduli space is really ill-behaved. E.g. the action of
$\OO(\Gamma)$ on $\kt_\Gamma^{\rm cpl}$ is not properly
discontinuous. Thus, $\kt_\Gamma^{\rm cpl}$ which already is not
Hausdorff, becomes even worse when divided out by $\OO(\Gamma)$
(cf.\ the discussion in Section \ref{Discrete}).

\bigskip

 There is yet another approach to these moduli
spaces where one actually fixes the underlying manifold and
constructs the moduli space as a quotient of the space of hyperk{\"a}hler
metrics by the diffeomorphism
group. We will briefly discuss this.

\medskip

Let $M$ be a compact oriented differentiable manifold of real dimension
$4n$ and let $q_M$ be the quadratic form on $H^2(M,\IZ)$ given by
$q_M(\alpha)=c_n\cdot\int_M\alpha^2\sqrt{\hat A(M)}$. We write
$\Gamma=(H^2(M,\IZ),q_M)$ and call this identification $\varphi_0$.

By ${\rm Diff}(M)$ we denote the group of orientation-preserving
diffeomorphisms of $M$. In fact, at least for $b_2\ne6$, the group
 ${\rm Diff}(M)$ is
the full diffeomorphism group of $M$, as any
orientation-reversing diffeomorphism $f$ would induce an isomorphism
$(H^2(M,\IZ),q_M)\cong (H^2(M,\IZ),-q_M)$ which is impossible
for $b_2(M)\ne6$. The set of all hyperk{\"a}hler metrics $g$ on $M$
is denoted by ${\rm Met}^{\rm HK}(M)$. Clearly, ${\rm Diff}(M)$ acts 
naturally on ${\rm Met}^{\rm HK}(M)$ by $(f,g)\mapsto f^*g$.

\begin{definition} The group
${\rm Diff}_{\rm o}(M)\subset{\rm Diff}(M)$ is the connected component
of ${\rm Diff}(M)$ containing the identity
${\rm id}_M\in{\rm Diff}(M)$.
The group
${\rm Diff}_{\rm *}(M)\subset{\rm Diff}(M))$ is the kernel of the natural
representation ${\rm Diff}(M)\to\OO(H^2(M,\IZ),q_M)$.
\end{definition}

Mapping $g\in{\rm Met}^{\rm HK}(M)$ to $(M,g,\varphi_0)\in\kt_\Gamma^{\rm
  met}$ induces a commutative diagram

$$\begin{array}{ccc}
{\rm Met}^{\rm HK}(M)/{\rm Diff}_*(M)&\rpfeil{5}{\eta}&\kt_\Gamma^{\rm
  met}\\
\downarrow&&\downarrow\\
{\rm Met}^{\rm HK}(M)/{\rm Diff}(M)&\rpfeil{5}{\bar\eta}&\km_\Gamma^{\rm met}\\
\end{array}$$

Note that $\eta$ is well-defined. Indeed, if $f\in{\rm Diff}_*(M)$,
then $(M,g,\varphi_0)\sim(M,f^*g,\varphi_0\circ f^*)=(M,f^*g,\varphi_0)$.

\begin{remark}
It seems essentially
nothing is known about the quotient of the natural inclusion
${\rm Diff}_{\rm o}(M)\subset{\rm Diff}_*(M)$, not even for K3 surfaces,
i.e.\ $n=1$.
\end{remark}

Clearly, the image of $\eta$ (and $\bar\eta$) can contain only those HK
$(M',g',\varphi)\in\kt_\Gamma^{\rm met}$ whose underlying real manifold
$M'$ is diffeomorphic to $M$. Let $\kt_\Gamma^{\rm met}(M)$ and
$\km_\Gamma^{\rm met}(M)$ denote the union of all those connected
components.

\medskip

{\bf i)} In general, $\eta:{\rm Met}^{\rm HK}(M)/{\rm Diff}_*(M)\to\kt^{\rm
met}_\Gamma(M)$ is injective, but not surjective.\\
The injectivity is clear. Let us explain why surjectivity
fails in general.
If $(M,g,\varphi)\in{\rm
Im}(\eta)$ and $\psi\in\OO(H^2(M,\IZ),q_M)$, then
$(M,g,\varphi_0\circ\psi)\in{\rm Im}(\eta)$ if and only if there exists
$f\in{\rm Diff}(M)$ with $f^*=\psi$ but ${\rm Diff}(M)\to
\OO(H^2(M,\IZ),q_M)$ is not necessarily surjective.
E.g.\ for K3 surfaces the image does not contain $-{\rm id}$ and, more
precisely,
$\OO(H^2(M,\IZ),\cup)/{\rm Diff}(M)\cong\IZ/2\IZ$ (cf.\ \ref{diffeo}).
However in this
case the situation is rather simple, as $\kt^{\rm met}_\Gamma$ consists
of two components, interchanged by ${\rm id}_{H^2}$, and ${\rm Met}^{\rm HK}(M)/{\rm Diff}_*(M)$ is one of them.
For higher dimensional HK nothing is known about the image of 
${\rm Diff}(M)\to\OO(H^2(M,\IZ),q_M)$.

\medskip

{\bf ii)} The map $\bar\eta:{\rm Met}^{\rm HK}(M)/{\rm
  Diff}(M)\to\km_\Gamma^{\rm met}(M)$ is bijective.\\
 Indeed, if
  $(M,g,\varphi)\in\kt_\Gamma^{\rm met}(M)$, then $[(M,g,\varphi_0)]=
[(M,g,(\varphi_0\varphi^{-1})\varphi)]=(\varphi_0\varphi^{-1})[(M,g,\varphi)]=[(M,g,\varphi)]\in\km_\Gamma^{\rm
  met}(M)$. Thus, $\bar\eta$ is surjective. If
  $\bar\eta(M,g)=\bar\eta(M,g')$, then there exists $\psi\in\OO(\Gamma)$
  such that $(M,g,\varphi_0)\sim(M,g',\psi\circ\varphi_0)$ and hence there
  exists $f\in{\rm Diff}(M)$ with $f^*g=g'$ (note that for $b_2(M)=6$ one
  would have to argue that $f$ can be chosen orientation-preserving)
and $\varphi_0=\psi\circ\varphi_0\circ f^*$. Thus, $[(M,g)]=[(M,g')]$ in
${\rm Met}^{\rm HK}(M)/{\rm Diff}(M)$ and hence $\bar\eta$ is injective.

\bigskip

One last word concerning the stabilizer of the action of ${\rm Diff}(M)$.
Clearly, the stabilizer of a hyperk{\"a}hler metric $g$ is the isometry group
${\rm Isom}(M,g)$ of $(M,g)$. This group is compact (cf.\ \cite{Besse}).
Hence the stabilizer of $g\in{\rm Met}^{\rm HK}(M)$ is a compact group.
Moreover, ${\rm Isom}(M,g)\cap{\rm Diff}_*(M)$ is finite. Indeed, if
$f\in{\rm Isom}(M,g)$ then $f$ maps any $g$-compatible complex structure $I$
to another $g$-compatible complex structure $f^*I$.
If in addition $f^*={\rm id}$ on $H^*(M,\IR)$, then the map $I\mapsto f^*I$ must also be the identity.
Hence $f\in{\rm Aut}(M,I)$ for any $g$-compatible complex structure $I$. Since
 $H^0((M,I),\kt)=0$, the latter group is discrete and, therefore,
${\rm Aut}(M,I)\cap{\rm Isom}(M,g)$
is finite. Hence, the action of ${\rm Diff}_*(M)$
on ${\rm Met}^{\rm HK}(M)$ has finite stabilizer.
 
\section{Period domains}\label{ps}

The moduli spaces that have been introduced in the last section
will be studied by means of various period maps. In this section
we define and discuss the spaces in which these maps take their
values, the period domains.

Let $\Gamma$ be a lattice of signature $(m,n)$. The standard
example for $\Gamma$ is the K3 lattice $2(-E_8)\oplus 3U$, where
$U$ denotes the hyperbolic plane $\left( \IZ^2 ,
\bigl(\begin{smallmatrix}0&1\\1&0\end{smallmatrix}\bigr)\right)$.
However, $\Gamma$ might in general be non-unimodular. This will be
of no importance in this section, as only the real vector space
$\Gamma_\IR:=\Gamma\otimes\IR$ is going to be used. In fact,
usually we will work with an arbitrary vector space $V$, but
$\Gamma$ will nevertheless occur in the notation. I hope this
will not lead to any confusion.
\subsection{Positive (oriented) subspaces}\label{posor}

 \noindent
Let $V$ be a real vector space  that is endowed with a bilinear
form $\langle \ ,\  \rangle$ of signature $(m,n)$, e.g.\
$V=\Gamma_\IR$. We will also write $x^2 $ for $\langle
x,x\rangle$. Fix $k\le m$ and consider the space of all
$k$-dimensional subspaces $W\subset V$ such that $\langle \ ,\
\rangle$ restricted to $W$ is positive definite. We will denote
this space by ${\Gr}^{\rm p}_k (V)$. Clearly, ${\Gr}^{\rm p}_k
(V)$ is an open non-empty subset of the Grassmannian ${\Gr}_k(V)$.

In order to describe ${\Gr}^{\rm p}_k (V)$ as a homogeneous space
we consider the natural action of
$\OO(V)$ on ${\Gr}_k^{\rm p}(V)$ given by $(\varphi ,W)\mapsto \varphi
(W)$. The stabilizer of a point $W_0 \in  \Gr^{\rm p}_k (V)$
is $\OO(W_0) \times \OO(W_0^\perp)$. Since the action is
transitive, one obtains the following description

$$\framebox{$\Gr^{\rm p}_k (V) \cong  \raisebox{.3ex}{$\OO(V)$}/\raisebox{-.6ex}{$\OO(W_0)
\times \OO(W_0^\perp )$}\cong
\raisebox{.3ex}{$\OO(m,n)$}/\raisebox{-.6ex}{$\OO(k) \times
\OO(m-k,n)$}$}$$

The second isomorphism depends on the choice of a basis of the
spaces $W_0$ and $W_0^\perp$.

Next consider the space $\Gr^{\rm po}_k(V)$ of all oriented positive
subspaces $W\subset V$ of dimension $k$. Clearly, the natural map
$\Gr_k^{\rm po}(V)\to\Gr_k^{\rm p}(V)$ is a $2:1$ cover.
 Again, $\OO(V)$ acts transitively on $\Gr_k^{\rm po} (V)$
and the stabilizer of an oriented positive subspace $W_0$ is
$\SO(W_0) \times \OO(W_0^\perp)$. Thus,

$$\framebox{$ \Gr_2^{\rm po} (V) \cong \raisebox{.3ex}{$\OO(V)$} /\raisebox{-.6ex}{$\SO(W_0) \times \OO(W_0^\perp
)$} \cong \,$\raisebox{.3ex}{$\OO(m,n)$}/\raisebox{-.6ex}{$\SO(k)
\times \OO(m-k,n)$}}$$

\subsection{Planes and complex lines}\label{planes}

For $k=2$ the space $\Gr_2^{\rm po}(V)$ allows an alternative description. It turns
out that  there is a natural bijection between this space and the space
$$Q_\Gamma:=\{x~|~x^2=0,\ (x+\bar x)^2>0\}\subset\IP(\Gamma_\IC),$$
where we use the $\IC$-linear extension of $\langle~,~\rangle$.
Note that the second condition in the definition of $Q_\Gamma$ is
well posed, i.e.\ independent of the representative
$x\in\Gamma_\IC$ of the line $x\in\IP(\Gamma_\IC)$, as long as
the first condition $x^2=0$ is satisfied. Clearly, $Q_\Gamma$ is
an open subset of a non-singular quadric hypersurface in
$\IP(\Gamma_\IC)$.

To any $x\in Q_\Gamma$ one associates the plane
$W_x:=\Gamma_\IR\cap(x\IC\oplus\bar x\IC)\subset\Gamma_\IR$
endowed with the orientation given by $({\rm Re}(x),{\rm
Im}(x))$. Since $x\IC\oplus\bar x\IC$ is invariant under
conjugation, this space is indeed a real plane. Moreover,
$W_{\lambda x}=\Gamma_\IR\cap(\lambda x\IC\oplus \bar\lambda\bar
x\IC)=W_x$ and $({\rm Re}(\lambda x),{\rm Im}(\lambda x))=({\rm
Re}(x),{\rm Im}(x))\bigl(\begin{smallmatrix}{\rm Re}(\Lambda)&{\rm
Im}(\lambda)&\\-{\rm Im}(\lambda)&{\rm
Re}(\Lambda)\end{smallmatrix}\bigr)$, where the matrix has
positive determinant. Hence, the oriented plane $W_x$ is
well-defined, i.e.\ it only depends on $x\in\IP(\Gamma_\IC)$. It
is positive, since $(\lambda x+\bar\lambda \bar
x)^2=\lambda\bar\lambda(x+\bar x)^2>0$ for $\lambda\ne0$.

 Conversely, if
$W\in\Gr_2^{\rm po} (\Gamma_\IR)$, then choose a positively
oriented orthonormal basis $w_1,w_2\in W$ and set $x:=w_1+iw_2$.
Then $W=W_x$ and $x^2=0$, $(x+\bar x)^2=(2w_1)^2>0$. Moreover,
$x\in\IP(\Gamma_\IC)$ does not depend on the choice of the basis
and any $x\in Q_\Gamma$ can be written in this form.

Thus, one has a bijection
$$\framebox{$Q_\Gamma\cong\Gr_2^{\rm po} (\Gamma_\IR)$}$$

\subsection{Planes and three-spaces}

For our purpose the spaces $\Gr^{\rm po}_2(\Gamma_\IR)$,
$\Gr^{\rm po}_3(\Gamma_\IR)$, and
$\Gr_4^{\rm po}(\Gamma_\IR\oplus U_\IR)$ are the most interesting ones. In the next
two sections we will study how they are related to each other. To this end
let us first introduce the space

$$\Gr_{2,1}^{\rm po}(\Gamma_\IR):= \{ (P,\omega )\
\vert \ P\in \Gr_2^{\rm po} (\Gamma_\IR )\ , \ \omega \in P^\perp
\subset \Gamma_\IR\ , \ \omega^2 >0 \}.$$

Clearly, this space projects naturally to $\Gr_2^{\rm po}(\Gamma_\IR)$
by  $(P,\omega ) \mapsto P$. The fibre over the point $P$ is the quadratic
cone $\{ \omega \ \vert \ \omega^2 > 0\} \subset P^\perp \subset \Gamma_\IR $.
If $\Gamma$ has signature $(3,b-3)$, this cone consists of exactly two
connected components, which can be identified with each other by $\omega\mapsto-\omega$.
Thus, the fibre of $\Gr_{2,1}^{\rm po}(\Gamma_\IR)\to\Gr_2^{\rm po}
(\Gamma_{\IR})$ over $P$ in this case is the disjoint union of two copies of a connected
cone, which will be called $\kc_P$.

In fact, $\Gr_{2,1}^{\rm po}(\Gamma_\IR)\to\Gr_2^{\rm po}(\Gamma_\IR)$
is a trivial cover, i.e.\ $\Gr_{2,1}^{\rm po}(\Gamma_\IR)$ splits into
two components. This can either be deduced from the fact that
$\Gr_{2}^{\rm po}(\Gamma_\IR)\cong\OO(3,b-3)/(\SO(2)\times\OO(1,b-3))$
is simply connected (cf.\ Section \ref{Top}) or from the following argument:
If we fix an oriented positive  three-space $F\in\Gr_3^{\rm po}(\Gamma_\IR)$,
then the orthogonal
projection $P\oplus\IR\omega\to F$ for any $\omega\in\pm\kc_X$ must be an
isomorphism, since $F^\perp$ is negative. Thus, we can single out
one of the two connected components of $\pm\kc_P$ by requiring that
$P\oplus\IR\omega\cong F$ is compatible with the orientations on both spaces.

 Mapping $(P,\omega )$ to the oriented positive 
three-space $F(P,\omega ):= P\oplus \omega\IR$ and the scalar
$\omega^2 \in \IR_{>0}$ defines a map $\Gr_{2,1}^{\rm po}(\Gamma_\IR)\to
\Gr^{\rm po}_3 (\Gamma_\IR ) \times \IR_{>0}$. The map
is surjective and the fibre over a point $(F,\lambda)$ can be
identified with the set of all $\omega \in F$ with $\omega^2 =
\lambda$ which is a two-dimensional sphere. \medbreak\noindent
Thus, one has the following diagram

$$\xymatrix{
 \Gr_{2,1}^{\rm po}(\Gamma_\IR)\ar[r]^{\scriptscriptstyle S^2~~~~~~~~~~~~~~~~~~~~~~~~~~~~~}\ar[d]^{\scriptscriptstyle \pm\kc_P}&\Gr_3^{\rm po} (\Gamma_\IR ) \times \IR_{>0}
\cong\left(\raisebox{.3ex}{$\OO(m,n)$} /
\raisebox{-.6ex}{$\SO (3)\times \OO(m-3,n)$}\right) \times \IR_{>0}&&\\
\Gr_2^{\rm
po}(\Gamma_\IR)\ar[r]^{\cong~~~~~~~~~~~~~~~~~~~~~~~~~~~~~~~}
&\raisebox{.3ex}{$\OO(m,n)$} / \raisebox{-.6ex}{$\SO (2)\times
\OO(m-2,n)$} ~~~~~~~~~~~~~~~~~~~~~~~~~~~~~~~~~~~&&}$$

 Note that the
two natural compositions $S^2 \subset\Gr_{2,1}^{\rm
po}(\Gamma_\IR) \to \Gr_2^{\rm po} (\Gamma_\IR)$ and
$\kc_P\subset\Gr_{2,1}^{\rm po}(\Gamma_\IR)\to\Gr_3^{\rm
p}(\Gamma_\IR) \times\IR_{>0}$ are both injective.


\subsection{Three- and four-spaces}\label{34spaces}

From now on we will assume that $\Gamma$ has signature $(3,b-3)$.
Furthermore, let us fix a standard basis $(w,w^*)$ of $U$, i.e.\
$w^2={w^*}^2=0$ and $\langle w,w^*\rangle=1$. We will see that
the space of four-spaces in $\Gamma_\IR\oplus U_\IR$ relates
naturally to the space of three-spaces in $\Gamma_\IR$.
Explicitly, we will show

$$\framebox{$\begin{array}{rcl}
\Gr^{\rm po}_3(\Gamma_\IR)\times\IR_{>0}\times\Gamma_\IR&\cong&
\Gr^{\rm po}_4(\Gamma_\IR\oplus U_\IR)\\
&\cong& \raisebox{.3ex}{$\OO(4,b-2)$} /
\raisebox{-.6ex}{$\SO (4)\times \OO(b-2)$}
\end{array}$}$$

The second isomorphism follows from Section \ref{posor}. The first one is
given as follows.

$$\phi:(F,\alpha,B)\mapsto\Pi:=B'\IR \oplus F' ,$$
where $F':=\{f-\langle f, B\rangle w|f\in F\}$ and
$B':=B+\frac{1}{2}(\alpha-B^2)w+w^*$. Clearly, $\langle f-\langle
f,B\rangle w, B'\rangle=\langle f-\langle f, B\rangle
w,B+w^*\rangle=0$ and thus the decomposition is orthogonal.
Furthermore, $(f-\langle f, B\rangle w)^2=f^2>0$ for $0\ne f\in
F$ and $B'^2=B^2+\alpha-B^2=\alpha>0$. Hence, $\Pi$ is a positive
four-space. Its orientation is induced by the orientation of $F\cong F'$ and
the decomposition $\Pi=B'\IR \oplus F'$.

In order to see that $\phi$ is bijective we study the inverse map
$\psi:\Pi\mapsto(F,B'^2,B)$, where $F$, $B'$, and $B$ are defined
as follows: One first introduces $F':=\Pi\cap w^\perp$. This
space is of dimension three, since otherwise $\Pi\subset
w^\perp=\Gamma_\IR\oplus w\IR$ and the latter space does not
contain any positive four-space. Again by the positivity of $\Pi$
one finds $w\not\in F'\subset \Pi$. Hence,
$F:=\pi(F')\subset\Gamma_\IR$ is a positive three-space, where
$\pi:\Gamma_\IR\oplus U_\IR\to\Gamma_\IR$ is the natural
projection. Furthermore, there exists a $B'\in\Pi$ such that
$\Pi=B'\IR \oplus F'$ is an orthogonal splitting. As before $B'$
cannot be contained in $w^\perp$. Thus, one can rescale $B'$ such
that $\langle B',w\rangle=1$. This determines $B'$ uniquely. Since
$B'\in\Pi$, one has $B'^2>0$. The B-field is by definition
$B:=\pi(B')$. One easily verifies that $\psi$ and $\phi$ are
indeed inverse to each other.

\subsection{Pairs of planes}\label{PP}

The last space we will discuss in this series of period domains is the space
of orthogonal oriented positive planes in $\Gamma_\IR\oplus U_\IR$,
i.e.

$$\Gr_{2,2}^{\rm po}(\Gamma_\IR\oplus U_\IR)=\{(H_1,H_2)\ |\ H_i\in\Gr_2^{\rm po}(\Gamma_\IR\oplus U_\IR), H_1\perp H_2\}.$$

Using the same techniques as before this space can also be
described as an homogeneous space as follows

$$\framebox{$\begin{array}{rcl}
 \Gr_{2,2}^{\rm po}(\Gamma_\IR\oplus
U_\IR)&\cong & \raisebox{.3ex}{$\OO(\Gamma_\IR\oplus U_\IR)$} /
\raisebox{-.6ex}{$\SO (H_1)\times \SO(H_2)\times \OO((H_1\oplus H_2)^\perp)$}\\
&\cong& \raisebox{.3ex}{$\OO(4,b-2)$} / \raisebox{-.6ex}{$\SO
(2)\times \SO(2)\times \OO(b-2)$},
\end{array}$}$$

for some chosen point
$(H_1,H_2)\in \Gr_{2,2}^{\rm po}(\Gamma_\IR\oplus U_\IR)$,
respectively basis of the spaces $H_1$, $H_2$, and $(H_1\oplus H_2)^\perp$.

We will be interested in the natural projection
$$\framebox{$\pi:\Gr_{2,2}^{\rm po}(\Gamma_\IR\oplus U_\IR)\twoheadrightarrow\Gr_4^{\rm po}(\Gamma_\IR\oplus U_\IR),~
(H_1,H_2)\mapsto\Pi:=H_1\oplus H_2
$}$$
and in the injection
$$\framebox{$\gamma:\Gr_{2,1}^{\rm po}(\Gamma_\IR)\times\Gamma_\IR\hookrightarrow \Gr_{2,2}^{\rm po}(\Gamma_\IR\oplus U_\IR)$}$$
which is compatible with
$\Gr_{2,1}^{\rm po}(\Gamma_\IR)\to\Gr_3^{\rm po}(\Gamma_\IR)$.

Let us first study the projection. Using the above description of
both spaces as homogeneous spaces this map corresponds to
dividing by $\SO(4)/(\SO(2)\times\SO(2)) $. The fibre of $\pi$
over $\Pi\in\Gr_4^{\rm po}(\Gamma_\IR\oplus U_\IR)$ is canonically
isomorphic to $\Gr_2^{\rm po}(\Pi)$ via $(H_1,H_2)\mapsto H_1$.
The inverse image of $H\in \Gr_2^{\rm po}(\Pi)$ is $(H,H^\perp)$,
where $H^\perp$ gets its orientation from $\Pi$ and the
decomposition $\Pi=H\oplus H^\perp$.

Thus, one obtains the following description of the fibre
$$\pi^{-1}(\Pi)\cong\Gr_2^{\rm po}(\Pi)\cong S^2\times S^2.$$
The second isomorphism is derived as in Section \ref{planes} from
$$\Gr_2^{\rm po}(\Pi)=\{x\in\IP(\Pi_\IC)\ |\ x^2=0\}\cong\IP^1\times\IP^1.$$
Note that $(x+\bar x)^2>0$ is automatically satisfied, for $\langle\, \, ,\,\, \rangle$ on $\Pi$ is positive  by assumption.

Let us now turn to the injection $\gamma$, which is defined as
follows. We set $\gamma((P,\omega),B)=(H_1,H_2)$ with
$$H_1:=\{x-\langle x, B\rangle w\ |\ x\in P\}$$
and
$$H_2:=\left(\frac{1}{2}(\alpha-B^2)w+w^*+B\right)\IR
\oplus\left(\omega-\langle\omega,B\rangle w\right)\IR,$$ where as
before $(w,w^*)$ is the standard basis of $U$  and
$\alpha=\omega^2$.

The isomorphism $P\cong H_1$, $x\mapsto x-\langle x, B\rangle w$
endows $H_1$ with an orientation. A natural orientation of $H_2$ is given by
definition. Observe that $H_1$ only depends on $P$ and $B$, whereas
$H_2$ depends on
$\omega$ and $B$. One easily verifies that the map $\gamma$ is injective and
that it commutes with the projections to
$$\Gr_3^{\rm po}(\Gamma_\IR\oplus U_\IR)\times\IR_{>0}\times\Gamma_\IR\cong\Gr_4^{\rm po}(\Gamma_\IR\oplus U_\IR).$$

Recall that the fibre of $\Gr_{2,1}^{\rm po}(\Gamma_\IR)\times
\Gamma_\IR\to\Gr_4^{\rm po}(\Gamma_\IR\oplus U_\IR)$ is $S^2$, whereas the
fibre of $\pi:\Gr_{2,2}^{\rm po}(\Gamma_\IR\oplus U_\IR)\to \Gr_4^{\rm
  po}(\Gamma_\IR\oplus U_\IR)$ is $S^2\times S^2$. It can be checked that
the embedding $\gamma$ does not identify the fibre $S^2$ neither with the diagonal
nor with one of the two factors.
In algebro-geometric terms $S^2\subset S^2\times S^2$ is a hyperplane section of $\IP^1\times\IP^1$
with respect to the Segre embedding in $\IP^3$.

\begin{remark} Note that the projection $\Gr_{2,1}^{\rm
po}(\Gamma_\IR)\times\Gamma_\IR\to \Gr_{2,1}^{\rm
po}(\Gamma_\IR)\to \Gr_2^{\rm po}(\Gamma_\IR)$ does not extend,
at least not canonically, to a map $\Gr_{2,2}^{\rm
po}(\Gamma_\IR\oplus U_\IR)\to \Gr_2^{\rm po}(\Gamma_\IR)$.
Geometrically this will be interpreted by the fact that not any point
in the $(2,2)$-CFT moduli space of K3 surfaces  canonically defines
a complex structure. More recently, it has become clear that generalized K3 surfaces, a notion that relies
on Hitchin's generalized Calabi-Yau structures
\cite{Hitchin}, might be useful to give a geometric interpretation
to every $N=(2,2)$-SCFT (see \cite{HuyK3}) 
\end{remark}

We summarize the discussion of this paragraph in the following
commutative diagram

$$\xymatrix{\Gr_{2,2}^{\rm po}(\Gamma_\IR\oplus
U_\IR)\ar@{>>}[r]^{\scriptscriptstyle S^2\times S^2}&\Gr^{\rm po}_4(\Gamma_\IR\oplus U_\IR)\\
\Gr^{\rm po}_{2,1}(\Gamma_\IR)\times\Gamma_\IR\ar @{^{(}->}
[u]\ar @{>>}[r]^{\scriptscriptstyle~~~
S^2~~~~~~~~}\ar@{>>}[d]^{\scriptscriptstyle \Gamma_\IR}&\Gr_3^{\rm po}(\Gamma_\IR)\times \IR_{>0}\times\Gamma_\IR\ar@{=}[u]\ar@{>>}[d]^{\scriptscriptstyle\Gamma_\IR}\\
\Gr_{2,1}^{\rm po}(\Gamma_\IR)\ar@{>>}[r]^{\scriptscriptstyle
S^2~~~~~~}\ar@{>>}[d]^{\scriptscriptstyle \pm\kc_P}&\Gr_3^{\rm po}(\Gamma_\IR)\times \IR_{>0}\\
\Gr_2^{\rm po}(\Gamma_\IR)& }$$

\subsection{Calculations in the Mukai lattice}

We shall indicate how the formulae change if we pass to the Mukai bilinear form. This will enable
us to make the description of the various period spaces and period maps compatible with conventions used 
elsewhere. We include this discussion for completeness, but it 
is not necessary for the understanding of the later sections.

Clearly the hyperbolic lattice $U$ with the standard basis $w,w^*$ is isomorphic to $-U$ via
$w\mapsto -w,w^*\mapsto w^*$. This extends to a lattice isomorphism
$$\eta:\Gamma\oplus U\cong \Gamma\oplus(-U)=:\tilde\Gamma.$$
For any oriented four-manifold $M$ underlying a K3 surfaces we can identify $H^*(M,\IZ)$
endowed with the standard intersection pairing with $\Gamma\oplus U$ such that
$w^*=1\in H^0(M,\IZ)$, $w=[{\rm pt}]\in H^4(M,\IZ)$, and $\Gamma\cong H^2(M,\IZ)$. Then $\tilde\Gamma$ is naturally isomorphic to
$H^*(M,\IZ)$ with the Mukai pairing $(\alpha_0+\alpha_2+\alpha_4,\beta_0+\beta_2+\beta_4)_{\tilde\Gamma}=
-\alpha_0\beta_4-\alpha_4\beta_0+\alpha_2\beta_2$, where $\alpha_i,\beta_i\in H^i(M,\IZ)$.

The identification of $\Gamma\oplus U$ and $\tilde\Gamma$ with the cohomology of a K3 surface induces a ring structure
on both lattices, i.e.\ in both cases we define $(\lambda w^*+x+\mu w)^2:=\lambda^2 w^*+2\lambda x+(2\lambda \mu+x^2)w$.
Note that $\eta$ does not respect these ring structures.

Using the ring structure on $\tilde\Gamma_\IR$ we can let act any element $B_0\in\Gamma_\IR$
on $\tilde\Gamma_\IR$ via its exponential $\exp(B_0) =w^*+B_0+(B_0^2/2)w$.

\begin{lemma}
For any $B_0\in\Gamma_\IR$ one has $\exp(B_0)\in\OO(\tilde\Gamma_\IR)$.
\end{lemma}

\prf
This results from the following straightforward calculation
\begin{eqnarray*}
&&\left(\exp(B_0)\cdot(\lambda w^*+x+\mu w)\right)_{\tilde\Gamma}^2\\
&=&\left(\lambda w^*+(\lambda B_0+x)+(\mu+\langle B_0,x\rangle+\lambda \frac{B_0^2}{2})w\right)_{\tilde\Gamma}^2\\
&=&x^2-2\lambda\mu=\left(\lambda w^*+x+\mu w\right)_{\tilde\Gamma}^2.
\end{eqnarray*}
\qed

\bigskip
Later we shall study the map $\varphi_{B_0}$ associated 
to any $B_0\in \Gamma_\IR$ (see Section \ref{geomsym}. By definition
$\varphi_{B_0}\in\OO(\Gamma_\IR\oplus U_\IR)$
acts on $\Gamma_\IR\oplus U_\IR$ by
\begin{eqnarray*}
w\mapsto w,~~w^*\mapsto B_0+w^*-\frac{B_0^2}{2}w,\\
x\mapsto x-\langle B_0,x\rangle w, ~~{\rm for}~x\in \Gamma_\IR.
\end{eqnarray*}

Let us compare $\exp(B_0)$ with  $\varphi_{B_0}$,

\begin{proposition}
Under the isomorphism $\eta:\Gamma_\IR\oplus U_\IR\cong\tilde\Gamma_\IR$ the automorphism
$\varphi_{B_0}$ corresponds to the action of $\exp(B_0)$, i.e.\
$\eta\circ\varphi_{B_0}=\exp(B_0)\circ\eta$.
\end{proposition}

\prf By definition, $\exp(B_0)$ acts by
\begin{eqnarray*}
w\mapsto w,~~w^*\mapsto (w^*+B_0+\frac{B_0^2}{2}w)\cdot w^*=B_0+w^*+\frac{B_0^2}{2}w,\\
x\mapsto (w^*+B_0+\frac{B_0^2}{2}w)\cdot x=x+\langle B_0,x\rangle w, ~~{\rm for}~x\in \Gamma_\IR,
\end{eqnarray*}
which yields the assertion.\qed

\bigskip

The isomorphism $\eta$ induces a natural isomorphism
$\Gr_{2,2}^{\rm po}(\Gamma_\IR\oplus U_\IR)\cong\Gr_{2,2}^{\rm po}(\tilde\Gamma_\IR)$. In order, to describe
the image $\eta(H_1,H_2)$ we will use the identification $Q_{\tilde\Gamma}\cong\Gr_2^{\rm po}(\tilde\Gamma_\IR)$
established in Section \ref{planes}. The positive plane associated to an element $x\in\tilde\Gamma_\IC$ with
$[x]\in Q_{\tilde\Gamma}$ will be denoted by $P_x$, i.e.\ $P_x$ is spanned by ${\rm Re}(x)$ and ${\rm Im}(x)$.
Clearly, $\exp(B)\cdot P_x=P_{\exp(B)x}$.

Let $(P,\omega)\in \Gr_{2,1}^{\rm po}(\Gamma_\IR)$ and $B\in\Gamma_\IR$. We denote 
$(H_1,H_2):=\gamma((P,\omega),0)$ and $(H^B_1,H^B_2):=\gamma((P,\omega),B)$.
Then a direct calculation shows $\varphi_B(H_1,H_2)=(H_1^B,H_2^B)$ and therefore
\begin{corollary}
$\eta(H_1^B,H_2^B)=\exp(B)\cdot\eta(H_1,H_2)=\exp(B)\cdot(P_\sigma,P_{\exp(i\omega)})=(P_{\exp(B)\sigma},P_{\exp(B+i\omega)})$.
\end{corollary}

\prf The only thing that needs a proof is $\eta(H_2)=P_{\exp(i\omega)}$. But this follow immediately from the
definition of $H_2$.\qed

\bigskip

In particular, via $\eta$ the image of
$\gamma:\Gr_{2,1}^{\rm po}(\Gamma_\IR)\times\Gamma_\IR\to\Gr_{2,2}^{\rm po}(\Gamma_\IR\oplus U_\IR)$ can be identified with 
$\exp(\Gamma_\IR)\cdot\left(\eta(\gamma(\Gr_{2,1}^{\rm po}(\Gamma_\IR)))\right)$.
\subsection{Topology of period domains}\label{Top}

Let us study some basic aspects of the topology of the period domains that
are of interest for us. Let $\Gamma$ be a lattice of signature
$(3,b-3)$.
We will consider the spaces:

$$\begin{array}{rcl}
\Gr_2^{\rm po}(\Gamma_\IR)&\cong&\raisebox{.3ex}{$\OO(3,b-3)$} /
\raisebox{-.6ex}{$\SO(2)\times \OO(1,b-3)$}\\
\Gr_4^{\rm po}(\Gamma_\IR\oplus U_\IR)&\cong&\raisebox{.3ex}{$\OO(4,b-2)$} / \raisebox{-.6ex}{$\SO
(4)\times \OO(b-2)$}\\
\Gr_3^{\rm po}(\Gamma_\IR)&\cong&\raisebox{.3ex}{$\OO(3,b-3)$} / \raisebox{-.6ex}{$\SO
(3)\times \OO(b-3)$}\\
\Gr_{2,2}^{\rm po}(\Gamma_\IR\oplus U_\IR)&\cong&\raisebox{.3ex}{$\OO(4,b-2)$} / \raisebox{-.6ex}{$\SO
(2)\times \SO(2)\times \OO(b-2)$}\\
\end{array}
$$

For simplicity we will suppose that $b>3$.

\begin{lemma}
The group $\OO(k,\ell)$ with $k,\ell>0$ has exactly four connected components.
\end{lemma}

\prf Write $\OO:=\OO(k,\ell)$. Then there are the following disjoint
unions $\OO=\OO^+\cup\OO^-$, $\OO=\OO_+\cup\OO_-$, and $\OO=\OO^+_+\cup
\OO^+_-\cup\OO^-_+\cup\OO^-_-$. Here, $\OO^\pm_\pm$ are defined as follows:
Write $\IR^{k+\ell}=W_0\oplus W_0^\perp$ with $W_0\subset \IR^{k+\ell}$ a
maximal positive subspace, which is endowed with an orientation. Then let $\OO^+$ and
$\OO^-$ (respectively, $\OO_+$ and $\OO_-$) be the subsets of all linear
maps $A\in\OO$ such that the orthogonal projection $AW_0\to W_0$
(respectively, $AW_0^\perp\to W_0^\perp$) is orientation preserving
resp.\ orientation reversing. By definition $\OO^+_+=\OO^+\cap \OO_+$, etc.
For any $A_0\in\OO^\pm_\pm$ the map $\OO^+_+\to\OO^\pm_\pm$,
$A\mapsto AA_0$ defines a homeomorphism. Thus, it suffices to show that
$\OO^+_+$ is connected.
\qed

\bigskip

Note that $\OO^+_+(k,l)$ is the
connected component of the identity.   It will thus also be denoted
$\OO_{\rm o}(m,n)$.

\begin{corollary} The space
$\Gr_2^{\rm po}(\Gamma_\IR)$ is connected, whereas the spaces
$\Gr_{2,1}^{\rm po}(\Gamma_\IR)$, $\Gr_3^{\rm po}(\Gamma_\IR)$, $\Gr_4^{\rm po}(\Gamma_\IR\oplus U_\IR)$, and $\Gr_{2,2}^{\rm po}(\Gamma_\IR\oplus U_\IR)$ consist of two connected components.
\end{corollary}

\prf Use the obvious fact that the inclusion
$\SO(2)\times\OO(1,b-3)\subset\OO(3,b-3)$ respects the decomposition
into connected components, i.e.\ $\OO^\pm_\pm(1,b-3)\subset\OO^\pm_\pm(3,b-3)$. Thus,
$\pi_0(\SO(2)\times\OO(1,b-3))\cong\pi_0(\OO(3,b-3))$. Similarly
for $\Gr_3^{\rm po}(\Gamma_\IR)$. Here $\OO(3,b-3)$ has four connected components,
but $\pi_0(\SO(3)\times\OO(b-3))=\IZ/2\IZ$, i.e.\ the components $\OO^-_\pm$ do not intersect the image of the inclusion. Hence,
$\pi_0(\Gr_3^{\rm po}(\Gamma_\IR))\cong\IZ/2\IZ$.
The remaining assertions are proved analogously.
\qed

\bigskip

We are also interested in the fundamental groups of these spaces.
In order to compute those, we recall the following classical facts.
\begin{proposition} One has
$\pi_1(\SO(2))=\IZ$, $\pi_1(\SO(k))=\IZ/2\IZ$ for $k>2$, and
$\pi_1(\OO_{\rm o}(k,\ell))\cong\pi_1(\SO(k))\times\pi_1(\SO(\ell))$.
\end{proposition}

\prf The first assertion follows from $\SO(2)\cong S^1$. The universal
cover of $\SO(k)$ for $k\geq3$ is the two-to-one cover
${\rm Spin}(k)\to\SO(k)$. The isomorphism in the last assertion is induced
by the natural inclusion
$\SO(k)\times\SO(\ell)\hookrightarrow\OO_{\rm o}(k,\ell)$.\qed

\bigskip

\begin{corollary}
All the Grassmanians
$\Gr_2^{\rm po}(\Gamma_\IR)$, $\Gr_{2,1}^{\rm po}(\Gamma_\IR)$, $\Gr_3^{\rm po}(\Gamma_\IR)$,
$\Gr_4^{\rm po}(\Gamma_\IR\oplus U_\IR)$, and \- $\Gr_{2,2}^{\rm po}(\Gamma_\IR\oplus U_\IR)$ are simply-connected, i.e.\ every connected component is simply
connected.
\end{corollary}

\prf
Since $$\Gr_2^{\rm po}(\Gamma_\IR)=\raisebox{.3ex}{$\OO(3,b-3)$} /
\raisebox{-.6ex}{$\SO(2)\times \OO(1,b-3)$}\cong\raisebox{.3ex}{$\OO_{\rm o}(3,b-3)$} /
\raisebox{-.6ex}{$\SO(2)\times \OO_{\rm o}(1,b-3)$},$$
we may use the exact sequence
$$\pi_1(\SO(2)\times\OO_{\rm o}(1,b-3))\rpfeil{5}{a}\pi_1(\OO_{\rm o}
(3,b-3))\to\pi_1(\Gr_2^{\rm po}(\Gamma_\IR))\to\pi_0(~~)\cong\pi_0(~~).$$
The map $a$ is compatible with the natural isomorphisms
$\pi_1(\SO(2)\times\OO_{\rm o}(1,b-3))\cong\pi_1(\SO(2))\times
\pi_1(\OO_{\rm o}(1,b-3))\cong\pi_1(\SO(2))\times\pi_1(\SO(1))\times\pi_1(\SO(b-3))$,
$\pi_1(\OO_{\rm o}(3,b-3))\cong\pi_1(\SO(3))\times\pi_1(\SO(b-3))$ and the natural
maps $\IZ\cong\pi_1(\SO(2))\times\pi_1(\SO(1))\to\pi_1(\SO(3))\cong\IZ/2\IZ$. Thus, $a$ is surjective and hence $\pi_1(\Gr_2^{\rm po}(\Gamma_\IR))=0$. The other assertions are proved analogously.
\qed

\bigskip

\begin{remark}
  Eventually, we list the real dimensions of our period spaces, which
  can easily be computed starting from
  $\Gr_2^{\rm po}(\Gamma_\IR)\cong Q_\Gamma$.
We have $\dim\Gr_2^{\rm po}(\Gamma_\IR)=2(b-2)$, $\dim\Gr_{2,1}^{\rm
  po}(\Gamma_\IR)=3(b-2)$, 
$\dim\Gr_3^{\rm po}(\Gamma_\IR)=3(b-3)$,
$\Gr_4^{\rm po}(\Gamma_\IR\oplus U_\IR)=4(b-2)$,
and $\dim\Gr_{2,2}^{\rm po}(\Gamma_\IR\oplus U_\IR)=4(b-1)$.
\end{remark}

\subsection{Density results}\label{dense}

Here we shall be interested in those points $P\in Q_\Gamma$ whose
orthogonal complement $P^\perp\subset\Gamma_\IR$ contains integral elements
$\alpha\in\Gamma$ of given length.  For simplicity we shall assume that
 $\Gamma$ is the K3 lattice $2(-E_8)\oplus 3U$, but all we will use is that
$\Gamma$ is even of index $(3,b-3)$ and that any primitive isotropic
element of $\Gamma$ can be complemented to a sublattice of $\Gamma$ which
is isomorphic to the hyperbolic plane. First note the following easy fact.

\begin{lemma}
If $0\ne\alpha\in\Gamma_\IR$ then $\alpha^\perp\cap Q_\Gamma$ is not empty.
\end{lemma}

\prf Indeed, $\alpha^\perp\subset\Gamma_\IR$ is a hyperplane containing at
least two linearly independent orthogonal positive vectors $x,y$. Thus, $P:=\langle
x,y\rangle\in\alpha^\perp\cap Q_\Gamma$.\qed

\bigskip

The quadric in $\IP(\Gamma_\IC)$ defined by the quadratic form
$\langle~,~\rangle$ on $\Gamma$ will be denoted $Z$, its real points form the
set $Z_\IR=\IP(\Gamma_\IR)\cap Z$.

\begin{proposition}\label{denseorb}
Let $0\ne\alpha\in\Gamma$. Then the set 
$$\bigcup_{g\in\OO(\Gamma)}g(\alpha^\perp\cap
Q_\Gamma)=\bigcup_{g\in\OO(\Gamma)}g(\alpha)^\perp\cap Q_\Gamma$$
is dense in $Q_\Gamma$.
\end{proposition}

\prf
We start out with the following observation: Let $\Gamma=\Gamma'\oplus U$
be an orthogonal decomposition and let $(v,v^*)$ be a standard basis of
$U$. For $B\in\Gamma'$ with $B^2\ne0$ we define $\varphi_B\in\OO(\Gamma)$
by $\varphi_B(v)=v$, $\varphi_B(v^*)=B+v^*-B^2/2\cdot v$, and
$\varphi_B(x)=x-\langle B,x\rangle v$ for $x\in\Gamma'$. It is easy to see
that indeed with this definition $\varphi_B\in\OO(\Gamma)$.
(We shall study a similarly defined automorphism
$\varphi_B\in\OO(\Gamma\oplus U)$ in Section \ref{Discrete}).

This automorphism has the remarkable property that for any $y\in\Gamma_\IR$
one has
$$\lim_{k\to\infty}\varphi^k_N[y]=[v]\in\IP(\Gamma_\IR).$$
In particular, we find that in the closure of the orbit
$\OO:=\OO(\Gamma)\cdot[\alpha]\subset\IP(\Gamma_\IR)$ there exists an isotropic
vector, i.e.\ $\overline\OO\cap Z_\IR\ne\emptyset$.

In order to prove the assertion of the proposition we have to show that for
any $P\in Q_\Gamma$ there exists an automorphism $g\in\OO(\Gamma)$ such that
 $g(\alpha)$ is arbitrarily close to $P^\perp$.
Indeed, in this case we find a codimension two subspace
$W\subset\Gamma_\IR$ close to $P^\perp$ containing $g(\alpha)$ and,
therefore,
 $W^\perp\in Q_\Gamma$ is close to $P$ and orthogonal to $g(\alpha)$.

Since $P^\perp$ contains some isotropic vector, it suffices to show that
any vector $[y]\in Z_\IR\subset \IP(\Gamma_\IR)$ is contained in $\overline\OO$.
As explained before, $\overline\OO\cap Z_\IR\ne\emptyset$. On the other hand,
$\overline\OO\cap Z_\IR$ is closed and $\OO(\Gamma)$-invariant. Thus, it
suffices to show that any $\OO(\Gamma)$-orbit $\OO_y:=\OO(\Gamma)\cdot[y]\subset
Z_\IR$ is dense. This is proved in two steps.

i) The closure $\overline\OO_y$ contains the subset $\{[x]\in
Z~|~x\in\Gamma\}$. Indeed, for any $x\in\Gamma$ primitive with $x^2=0$ one
finds an orthogonal decomposition $\Gamma=\Gamma'\oplus U$ with $x=v$,
where $(v,v^*)$ is a standard basis of the hyperbolic plane $U$. If we
choose
$B\in\Gamma'$ with $B^2\ne0$, then $\lim_{k\to\infty}\varphi_B^k[y]=[v]=[x]$, as we have
seen before. Hence, $[x]\in\overline\OO_y$.

ii) The set $\{[x]\in Z~|~x\in\Gamma\}$ is dense in $Z$. Indeed, if we
write $\Gamma=\Gamma'\oplus U$ as before, then the dense open subset
$V\subset Z_\IR$ of points of the form $[x'+\lambda v+v^*]$ with
$\lambda\in\IR$,
$x'\in\Gamma'_\IR$ is the affine quadric
$\{(x',\lambda)~|~2\lambda+x'^2=0\}\subset\Gamma_\IR\times\IR$ and thus is
given as the graph of the rational polynomial $\Gamma'_\IR\to\IR$,
 $x'\mapsto-x'^2/2$. Therefore, the rational points are dense in $V$. 

Combining both steps yields the assertion.\qed
\bigskip

\begin{corollary}
For any $m\in \IZ$ the subset
$$\{P\in Q_\Gamma~|~{\rm there~exists~a~primitive}~\alpha\in\Gamma\cap
P^\perp~{\rm with~}\alpha^2=2m\}$$
is dense in $Q_\Gamma$.
\end{corollary}

\prf In order to apply the proposition we only have to ensure that there is
a primitive element $0\ne\alpha\in\Gamma$ with $\alpha^2=2m$.
If $(w,w^*)$ is the standard base of a copy of the
hyperbolic plane  $U$ contained in $\Gamma$, we can choose
$\alpha=w+mw^*$.\qed

\bigskip

In fact, if $\alpha_1,\alpha_2\in\Gamma$ are primitive elements with
$\alpha_1^2=\alpha_2^2$ then there exists an automorphism
$\varphi\in\OO(\Gamma)$ with $\varphi(\alpha_1)=\alpha_2$ (cf.\ \cite[Thm.2.4]{LP} or Remark \ref{LPrem}). Thus, the assertion of the corollary is essentially
equivalent to the proposition (see \cite{Periodes} page 111).
Note that for general HKs we don't know which values of $2m$ can be realized.

As a further trivial consequence, one finds that the set of those $P\in
Q_\Gamma$ such that $P^\perp\cap \Gamma\ne0$ is dense in $Q_\Gamma$. One
can now go on and ask for those $P\in Q_\Gamma$ such that $P^\perp\cap
\Gamma$ has higher rank. Those with maximal rank, i.e.\
${\rm rk}(P^\perp\cap\Gamma)={\rm rk}(\Gamma)-2$, are called exceptional. An
equivalent definition is

\begin{definition}
A period point $P\in Q_\Gamma$ is exceptional if $P\subset\Gamma_\IR$ is
defined over $\IQ$, i.e.\ $P\in Q_\Gamma\cap\IP(\Gamma_{\IQ(i)})$.
\end{definition}

Clearly, $P$ is exceptional if there exist linearly independent elements
$\alpha_1,\ldots,\alpha_{{\rm rk}(\Gamma)-2}\in\Gamma$ such that
$P\subset\alpha_i^\perp$ for all $i$.
Note that if $P\in Q_\Gamma$ is exceptional, the orthogonal complement
$P^\perp$ always contains a lattice vector $x\in\Gamma$ with $x^2>0$ (use
that $\Gamma$ has signature $(3,b-3)$).

Next we will prove that also the exceptional points are dense in
$Q_\Gamma$. For K3 surfaces one can add further restrictions.

\begin{definition} Let $\Gamma$ be the K3 lattice.
A period point $P\in Q_\Gamma$ is called \emph{exceptional Kummer} if
$P\subset\Gamma_\IR$ is defined over $\IQ$ and for all $x\in P\cap\Gamma$
one has $x^2\equiv 0\mod 4$. 
\end{definition}

\begin{proposition}
Let $\Gamma$ be the K3 lattice. Then the set of exceptional Kummer points
$P\in Q_\Gamma$ is a dense subset of $Q_\Gamma$.
\end{proposition}

\prf We first prove the following statement. Let $L$ be an arbitrary
lattice. Then the set
$$\{[x]~|~x\in L~{\rm
is~  primitive~and}~x^2\equiv0\mod4\}\subset\IP(L_\IR)$$ is empty or
dense. Indeed, if $[x]$ is contained in this set and $y\in L$ is arbitrary,
then $[x+N\cdot y]\in\IP(L_\IR)$ converges towards $[y]$ for
$N\to\infty$. Moreover,
$(x+N\cdot y)^2\equiv x^2\equiv 0\mod4$ if $N$ is even. If $y\in L$ is
primitive and $y\ne x$ then there exist arbitrarily large even $N$ such that
$x+N\cdot y$ is again primitive. Since the set of all $[y]$ with $y\in L$
primitive is dense in $\IP(L_\IR)$, this proves the assertion.

Now let $P\in Q_\Gamma$ be spanned by orthogonal vectors
$y_1,y_2\in\Gamma_\IR$. Then by what was explained before we can find
$x_1\in\Gamma$ primitive with $x^2_1\equiv0\mod 4$ such that $[x_1]$ is
arbitrarily close to $[y_1]\in\IP(\Gamma_\IR)$. Furthermore, choose $x_2\in
x_1^\perp\subset\Gamma$ primitive and arbitrarily close to $y_2\in
y_1^\perp$
with $x_2^2\equiv0\mod4$ and set $P':=(\IZ x_1\oplus \IZ x_2)_\IR$. 
Such an element $x_2$ can be found, as $x_1^\perp\subset \Gamma$ contains a
copy
 of the hyperbolic plane $U$ and thus an element whose square is divisible
 by four, e.g.\ $2v+v^*$, where $(v,v^*)$ is a standard basis of $U$.
Then
$P'$ is close to $P$ and $(ax_1+bx_2)^2=a^2x_1^2+b^2x_2^2\equiv0\mod4$.\qed

\bigskip

We leave it to the reader to modify the above proof to obtain

\begin{corollary}
Let $\Gamma$ be an arbitrary lattice of signature $(3,b-3)$. Then the set
of exceptional period points is dense in $Q_\Gamma$.\qed
\end{corollary}

\section{Period maps}\label{periodmaps}

The aim of this section is to compare the various moduli spaces
introduced in Section \ref{ms} with the period domains of Section
\ref{ps} via period maps $\kp^{\rm cpl}$, $\kp$, $\kp^{\rm met}$,
$\kp^{(2,2)}$, and $\kp^{(4,4)}$.
\subsection{Definition of the period maps}\label{defpmaps}
The period maps we are about to define will fit into the
following two commutative diagrams:

$$\xymatrix{\kp^{\rm cpl}:&\kt_\Gamma^{\rm
cpl} \ar[r]&\Gr_2^{\rm
po}(\Gamma_\IR)\cong Q_\Gamma~~\\
\kp:&\kt_\Gamma\ar@{>>}[u]\ar@{>>}[d]^{\scriptscriptstyle S^2}\ar[r]&\Gr_{2,1}^{\rm po}(\Gamma_\IR)\ar@{>>}[u]\ar@{>>}[d]^{{\scriptscriptstyle S^2}}~~~\\
\kp^{\rm met}:&\kt^{\rm met}_\Gamma\ar[r]&\Gr_3^{\rm
po}(\Gamma_\IR)\times\IR_{>0}}$$
and
$$\xymatrix{ \kp^{(2,2)}:&\kt^{(2,2)}_\Gamma
\ar@{>>}[d]^{{\scriptscriptstyle S^2}} \ar[r]^{}& \Gr_{2,1}^{\rm
po}(\Gamma_\IR)\times\Gamma_\IR \ar@{>>}[d]^{{\scriptscriptstyle
S^2}} \ar@{^{(}->}[r] &\Gr_{2,2}^{\rm
po}(\Gamma_\IR\oplus U_\IR)     \ar@{>>}[dl]^{\scriptscriptstyle S^2\times S^2}   \\
\kp^{(4,4)}:&\kt_\Gamma^{(4,4)} \ar[r]^{}& \Gr^{\rm
po}_4(\Gamma_\IR\oplus U_\IR)   }   $$

The latter should be compatible with the two diagrams
$$\xymatrix{\kt_\Gamma^{(2,2)}\ar@{>>}[d]\ar@{>>}[r]&\kt_\Gamma\ar@{>>}[d]\\
\kt_\Gamma^{(4,4)}\ar@{>>}[r]&\kt_\Gamma^{\rm met}}
~~~~~~\xymatrix{\Gr_{2,1}^{\rm
po}(\Gamma_\IR)\times\Gamma_\IR\ar@{>>}[d]\ar@{>>}[r]&\Gr_{2,1}^{\rm
po}(\Gamma_\IR)\ar@{>>}[d]\\
\Gr_4^{\rm po}(\Gamma_\IR\oplus U_\IR)\ar@{>>}[r]&\Gr_3^{\rm
po}(\Gamma_\IR)\times\IR_{>0}}$$ and the period maps $\kp$ and $\kp^{\rm met}$.

The definition of the maps $\kp$, $\kp^{\rm met}$, $\kp^{\rm
cpl}$, $\kp^{(2,2)}$, and $\kp^{(4,4)}$ is straightforward. Let
$(X,\alpha,\varphi)=(M,g,I,\varphi)\in\kt_\Gamma$ and $B\in
H^2(X,\IR)=H^2(M,\IR)$ a B-field. By $\sigma$ we denote a
generator of $H^{2,0}(X)$.

Then we set:
\begin{eqnarray*}
\kp^{\rm cpl}(X,\varphi)&=&[\varphi(\sigma)]\in
Q_\Gamma\subset\IP(\Gamma_\IC)\\
&=&\varphi\langle\re(\sigma),\im(\sigma)\rangle\in\Gr_2^{\rm
po}(\Gamma_\IR)\\
\kp(X,\alpha,\varphi) &=&\left(\kp^{\rm
cpl}(X,\varphi),\varphi(\alpha)\right)\in\Gr_{2,1}^{\rm
po}(\Gamma_\IR)\\
\kp^{\rm
met}(M,g,\varphi)&=&\left(\varphi(H^2_+(M,g)),q(M,g)\right)\in\Gr_3^{\rm
po}(\Gamma_\IR)\times\IR_{>0}\\
\kp^{(2,2)}(M,g,I,B,\varphi)&=&\left(\kp(M,g,I,\varphi),\varphi(B)\right)\in\Gr_{2,1}^{\rm
po}(\Gamma_\IR)\times\Gamma_\IR\\
\kp^{(4,4)}(M,g,B,\varphi)&=&\left(\kp^{\rm
met}(M,g,\varphi),\varphi(B)\right)\in\Gr_3^{\rm
po}(\Gamma_\IR)\times\IR_{>0}\times\Gamma_\IR\cong\Gr_4^{\rm
po}(\Gamma_\IR\oplus U_\IR)\\
\end{eqnarray*}

We leave it to the reader to verify that all period maps are
$\OO(\Gamma)$-equivariant and that one indeed obtains the
above commutative diagrams.

Also note that there is a natural $\OO(\Gamma\oplus U)$-action on the two
period domains $\Gr_{2,2}^{\rm po}(\Gamma_\IR\oplus U_\IR)$ and
$\Gr_4^{\rm po}(\Gamma_\IR\oplus U_\IR)$, but the image
of $\kp^{(2,2)}$ (or its closure) is not left invariant under this action.

\subsection{Geometry and period maps}

Without going too much into the details we collect in the
following some important results about period maps. In particular, we
will translate geometric results, like the Global Torelli Theorem into
global properties of the period maps.

\bigskip

{\bf Local Torelli.} {\it The map $\kp^{\rm cpl}:\kt_\Gamma^{\rm
cpl}\to Q_\Gamma$ is holomorphic and locally (in $\kt_\Gamma^{\rm
cpl}$) an isomorphism} (cf.\  \cite{Beauv}).

\medskip

Recall that $\kt^{\rm cpl}_\Gamma$ has a natural complex
structure, but that the underlying topological space is not
Hausdorff. On the other hand, $Q_\Gamma$ is an open subset of a
non-singular quadric in $\IP(\Gamma_\IC)$ and, therefore, a nice
complex manifold.

Of course, the Local Torelli Theorem in the above version
immediately carries over to the other period maps. Thus, $\kp$,
$\kp^{\rm met}$, $\kp^{(2,2)}$, and $\kp^{(4,4)}$ are all locally
injective. Since the Teichm{\"u}ller spaces $\kt_\Gamma$,
$\kt_\Gamma^{\rm met}$ $\kt_\Gamma^{(2,2)}$, and $\kt^{(4,4)}$
are all Hausdorff, this shows that except $\kp^{\rm cpl}$ all
period maps define covering maps on their open images.

\bigskip

{\bf Twistor lines.} {\it Under the period map $\kp^{\rm cpl}$ the twistor line
$\IP^1=c(m^{-1}(M,g,\varphi))\subset\kt_\Gamma^{\rm cpl}$ (cf.\ Proposition
\ref{twistorlines})
is identified with a quadric in some linear subspace $\IP^2\subset\IP(\Gamma_\IC)$.}

\medskip
 Indeed, the $\IP^2$ is given as
$\IP(\varphi(H^2_+(M,g)_\IC))\subset\IP(\Gamma_\IC)$.

\bigskip

{\bf Surjectivity of the period map.} {\it The map $\kp^{\rm
cpl}:\kt_\Gamma^{\rm cpl}\to Q_\Gamma$ maps every connected
component of $\kt^{\rm cpl}_\Gamma$ onto $Q_\Gamma$} (cf.\ 
\cite{Huy2}).

\medskip

Analogous statements for the other period maps do not hold. In
these cases the assertion has to be modified. To see this let us
look at the fibres of $\kt_\Gamma\to\kt_\Gamma^{\rm cpl}$ over
$(X,\varphi)$. By definition of $\kt_\Gamma$ this is the K{\"a}hler
cone $\kk_X$ which, via the period map $\kp$, is identified with an
open subcone of the positive cone $\kc_{\kp^{\rm
cpl}(X,\varphi)}$ which is just one of the two connected components of the
fibre of $\Gr_{2,1}^{\rm
po}(\Gamma_\IR)\to Q_\Gamma$ over $\kp^{\rm cpl}(X,\varphi)$. For
a very general marked IHS $(X,\varphi)\in\kt_\Gamma^{\rm cpl}$
the K{\"a}hler cone $\kk_X$ is maximal, i.e.\ $\kk_X=\kc_X$. Thus,
for those points $\kp$ maps the fibre of
$\kt_\Gamma\to\kt_\Gamma^{\rm cpl}$ bijectively onto one of the connected
components $\kc_P$ or $-\kc_P$ of the fibre of
$\Gr_{2,1}^{\rm po}(\Gamma_\IR)\to Q_\Gamma$ over $P=\kp^{\rm
cpl}(X,\varphi)$. For special marked IHS $(X,\varphi)$, which
usually (e.g.\ for K3 surfaces) nevertheless form a dense subset
of $\kt_\Gamma^{\rm cpl}$, the K{\"a}hler cone is strictly smaller.

\bigskip

{\bf Density of the image.} {\it The image of every connected
component of $\kt_\Gamma$ under the period map $\kp$ is dense in
the connected component of the
period domain $\Gr_{2,1}^{\rm po}(\Gamma_\IR)$ containing it. Analogous
statements hold true for $\kp^{\rm met}$, $\kp^{(2,2)}$, and
$\kp^{(4,4)}$.}

\medskip

Let us say a few words about how the density is proved and how the boundary
$\Gr_{2,1}^{\rm po}(\Gamma_\IR)\setminus\kp(\kt_\Gamma)$ can be
interpreted.

Since $\kp^{\rm cpl}$ is surjective, we may consider
$(X,\varphi)\in\kt^{\rm cpl}$ and study the fibre of $\Gr_{2,1}^{\rm
  po}(\Gamma_\IR)\to\Gr_2^{\rm po}(\Gamma_\IR)\cong Q_\Gamma$ over
$\kp^{\rm cpl}(X,\varphi)$, which is $\pm\varphi(\kc_X)$.
The $\pm$-sign distinguishes the two connected components
of $\Gr_{2,1}^{\rm po}(\Gamma_\IR)$. The
image of the fibre $\kt_\Gamma\to\kt_\Gamma^{\rm cpl}$ over $(X,\varphi)$
is the open subcone $\varphi(\kk_X)\subset\varphi(\kc_X)$. We will
discuss its boundary and its complement:
If $\alpha\in\kc_X$ is general, then there exists
$(X',\varphi')\in\kt_\Gamma^{\rm cpl}$
which cannot be separated from $(X,\varphi)$
such that $\kp(X,\varphi)=\kp(X',\varphi')$ and
$\varphi(\alpha)\in\varphi'(\kk_{X'})$ (see \cite{GHJ}).
(Moreover, $X$ and $X'$ are birational.) Thus, the disjoint union
$\bigcup\varphi(\kk_X)$ over all $(X,\varphi)$ in the same connected
component and with the same period $\kp(X,\varphi)\in Q_\Gamma$ is dense in 
$\varphi(\kc_X)$.

\medskip

For a point $\alpha\in\partial\varphi(\kk_X)$ there always exists a rational
curve $C\subset X$ with $\int_C\alpha=0$ (see
\cite{Boucksom}), i.e.\ under the degenerate
K{\"a}hler structure $\alpha$ the volume of the rational curve $C$ shrinks to
zero. Thus, points in the boundary of $\kp(\kk_X)$ should be thought of as 
singular IHS/HK which are obtained by contracting certain rational
curves. Unfortunately, neither are we able to make this statement more
precise nor do we know that any point $\alpha\in\varphi(\kc_X)$ is actually
contained in the closure of some $\varphi'(\kk_{X'})$, where
$(X',\varphi')$
is as above. However, for K3 surfaces the situation is much better
understood (cf.\ \cite{Kob}).

\bigskip

{\bf Projective IHS.} {\it
The set of projective marked IHS forms a countable dense
union of hyperplane section of $Q_\Gamma$. If $\km_\Gamma^{\rm proj}\subset\km_\Gamma$ denotes the set of all K{\"a}hler IHS for which the underlying IHS
is projective,
then $\km_\Gamma^{\rm proj}\to\km_\Gamma^{\rm met}$ is surjective.}

\medskip

In fact, due to a general projectivity criterion for surfaces and
an analogous result for IHS (cf.\ \cite{GHJ}) one 
knows that an IHS $X$ is projective if and only if there exists an integral
$(1,1)$-class $\alpha$ with $q(\alpha)>0$. Thus,
$(X,\varphi)\in\kt_\Gamma^{\rm cpl}$ is projective if and only if $\kp^{\rm cpl}(X,\varphi)$ is contained in a hyperplane
orthogonal to some $\alpha\in\Gamma$ with $\alpha^2>0$. As we have seen before,
the set of such periods is dense in the period domain $Q_\Gamma$.
Since the fibre of $\kt_\Gamma\to\kt^{\rm met}_\Gamma$ is identified with
a quadric curve $\IP^1\subset\IP(\Gamma_\IC)$ under the projection $\kt_\Gamma\to\kt_\Gamma^{\rm cpl}$
and as such is intersected non-trivially
by every such hyperplane, the fibre contains at least one K{\"a}hler $(X,\alpha,\varphi)$
with $X$ projective. In other words, for any hyperk{\"a}hler metric $g$
on a manifold $M$ at least one of the complex structures
$\lambda=aI+bJ+cK$ defines
a projective IHS. In fact, the set of projective IHS is also dense among
the $(M,\lambda)$.

\bigskip

{\bf  Finiteness.} {\it The induced period maps
$$\begin{array}{cccl}
\overline\kp:&\km_\Gamma&\to&\raisebox{-.6ex}{$\OO(\Gamma)$}\setminus\raisebox{.3ex}{$\Gr_{2,1}^{\rm po}(\Gamma_\IR)$}\\
&&&\\
\overline\kp^{\rm met}:&\km_\Gamma^{\rm
  met}&\to&\raisebox{-.6ex}{$\OO(\Gamma)$}\setminus\raisebox{.3ex}{$\Gr_3^{\rm
    po}(\Gamma_\IR)$}\times
\IR_{>0}\cong\raisebox{-.6ex}{$\OO(\Gamma)$}\setminus\raisebox{.3ex}{$\OO(3,b-3)$}/\raisebox{-.6ex}{$\SO(3)\times\OO(b-3)$}\times\IR_{>0}
\end{array}$$
 are finite trivial covers of their images, i.e.\ every
moduli space has only finitely many connected components and each connected
component is mapped bijectively onto its image.

The same holds for the period map
$$\overline\kp^{\rm cpl}:\km_\Gamma^{\rm
  cpl}\to\raisebox{-.6ex}{$\OO(\Gamma)$}\setminus
\raisebox{.3ex}{$Q_\Gamma$}\cong\raisebox{-.6ex}{$\OO(\Gamma)$}
\setminus\raisebox{.3ex}{$\OO(3,b-3)$}/\raisebox{-.6ex}{$\SO(2)\times\OO(1,b-3)$}$$
except for non-Hausdorff points in the fibers.}
\medskip

Note that e.g.\ $\kt_\Gamma^{\rm cpl}$ might {\it a priori}
have infinitely many
components. That this is no longer possible for the quotient
$\km_\Gamma^{\rm cpl}=\OO(\Gamma)\setminus\kt_\Gamma^{\rm cpl}$ is a
consequence of the finiteness result in \cite[Thm. 4.3]{ Huybrechtsfinite}
which says that there are only finitely many different deformation types of
IHS with the same BB--form $q_X$. 
Since $Q_\Gamma$ is simply connected and $\kp^{\rm cpl}$ is
surjective, the cover $\kp^{\rm cpl}$ has to be trivial. In fact, in order to
make this precise one first should construct the `Hausdorff reduction'
of $\kt^{\rm cpl}_\Gamma$ by identifying all points that cannot be separated from each other. This Hausdorff space then is an honest {\'e}tale cover of the simply connected space $Q_\Gamma$ and, therefore, consists of several copies
of $Q_\Gamma$.

\bigskip

We leave it to the reader to deduce similar statements for the maps
$\kp^{(2,2)}$ and $\kp^{(4,4)}$.

\begin{remark}
This is essentially all that is known in the general case. For {\bf
K3 surfaces}
however the above results can be strengthened considerably as follows.
The Global Torelli for K3 surfaces shows that $\kt_\Gamma^{\rm cpl}$
consists of two connected components which are identified with each other by
$(X,\varphi)\mapsto (X,-\varphi)$ and which are not distinguished by
$\kp^{\rm cpl}$. The two components are separated by the
map $\kp:\kt_\Gamma\to\Gr_{2,1}^{\rm po}(\Gamma_\IR)$, which is injective
in the case of K3 surfaces. Analogously, $\kp^{\rm met}$, $\kp^{(2,2)}$,
and $\kp^{(4,4)}$ are all injective.
\end{remark}

The density results of Section \ref{dense} together with the description of the
periods of our list of examples of K3 surfaces in Section \ref{basics}
and the above information
about the period maps (i.e.\ the Global Torelli Theorem) yield:

\begin{proposition}
The following three sets are dense in the moduli space of marked K3 
surfaces:
{\rm i)} $\{(X,\varphi)~|~X\subset\IP^3~{\rm is
  ~a~quartic~hypersurface}\}$,
  
{\rm ii)} $ \{(X,\varphi)~|~X~{\rm is~an~elliptic~K3~surface}\}$, {\rm and}

{\rm iii)} $\{(X,\varphi)~|~X~{\rm is~a(n~exceptional)~Kummer~surface}\}.$
\qed\end{proposition}

\subsection{The diffeomorphism group of a K3 surface}\label{diffeo}

\begin{proposition}
Let $X$ be a K3 surface. The image of the natural map $\rho:{\rm
  Diff}(X)\to\OO(H^2(X,\IZ),\cup)$ is the subgroup
$\OO^+(H^2(X,\IZ),\cup)$, which  is of index two.
\end{proposition}

Recall (cf.\ Section \ref{Top}) that $\OO^+$ is the group of all
$A\in\OO$ that preserve the orientation of positive three-space (but
not necessarily of a negative
$19$-space). The proposition is due to Borcea \cite{Borcea}, who showed
the inclusion $\OO^+\subset{\rm Im}(\rho)$, and Donaldson \cite{Don}, who
showed equality.
We only reproduce Borcea's argument here.

\prf First note the following. If $(X_t,\varphi_t)$ is a connected path in
$\kt_\Gamma^{\rm cpl}$, then there exists a sequence of diffeomorphisms $f_t:X_0\cong X_t$ such that $\varphi_0\circ f^*=\varphi_t$.

Let now $\varphi$ be any marking of $X$ and consider $(X,\varphi)\in\kt_\Gamma^{\rm cpl}$. By $\kt_0$ we denote the connected component of $\kt_\Gamma^{\rm cpl}$ that contains this point. Pick $A\in\OO^+(H^2(X,\IZ),\cup)$. Then
$A$ acts on $\kt_\Gamma^{\rm cpl}$ and $Q_\Gamma$ by $\varphi A\varphi^{-1}$
and the period map $\kp^{\rm cpl}:\kt_\Gamma^{\rm cpl}\to Q_\Gamma$ is
equivariant. Since the restriction of the period map
$\kp^{\rm cpl}$ yields a surjective map $\kt_0\to Q_\Gamma$, there exists
a marked K3 surface $(X',\varphi')$ with $\kp^{\rm cpl}(X',\varphi')=
A\kp^{\rm cpl}(X,\varphi)=\kp^{\rm cpl}(X,A\varphi)$.
If $X$ is a general K3 surface such that $\kk_X\cong\kc_X$,
then $\pm{\varphi'}^{-1}\circ
(\varphi A):H^2(X,\IZ)\cong H^2(X',\IZ)$ is an isomorphism of periods mapping
$\kk_X$ to $\kk_{X'}$. By the Global Torelli Theorem there exists a (unique)
isomorphism $g:X'\cong X$ such that $g^*=\pm{\varphi'}^{-1}\circ(\varphi A)$.
By the remark above we also find a diffeomorphism $f:X\cong X'$ such
that $\varphi\circ f^*=\varphi'$. Hence, $\varphi\circ f^* g^*=\pm(\varphi A)$ and thus $(g\circ f)^*=\pm A$ is realized by a diffeomorphism of $X$.
In fact, the sign must be ``$+$'', as $g^*$, $f^*$, and $A$ preserve the orientation of a positive three-space.

It remains to show that $-{\rm id}$ is not contained in the image and this
was done by Donaldson using zero-dimensional moduli spaces of stable bundles 
on a double cover of the projective plane.  \qed

\bigskip

\begin{remark}
In the proof above we used the assumption that $n=1$ twice: When we applied
the Global Torelli Theorem and, of course, when using Donaldson invariants.
The surjectivity which is also crucial holds true also for $n>1$. Somehow,
the use of the Global Torelli Theorem seems a little strong, as we have no
need to know that $g^*$ is induced by a biholomorphic map, a diffeomorphism
would be enough.
\end{remark}

In \cite{Namikawa} Namikawa constructs an example of two four-dimensional IHS $X$ and $X'$ together with an isomorphism of their periods which preserves
the K{\"a}hler cone, but such that $X$ and $X'$ are not even birational.
To be more precise, he considers generalized Kummer varieties $X={\rm K}_2(T)$ and $X'={\rm K}_2(T^*)$ 
 associated to a complex torus $T$ and its dual
$T^*$. As the moduli space of complex tori is connected, one can endow
$X$ and $X'$ with markings $\varphi$ respectively $\varphi'$ such that
$(X,\varphi)$ and $(X',\varphi')$ are contained in the same connected
component $\kt_0$ of $\kt_\Gamma^{\rm cpl}$. His example shows that
$\OO^+(\Gamma)$ does not preserve $\kt_0$, i.e.\ there exists $A\in\OO^+$
such that $(X',A\varphi')\not\in\kt_0$ (with
$\kp(X',A\varphi')=\kp(X,\varphi)$).
Indeed, after identifying non-separated points in $\kt_\Gamma^{\rm cpl}$
the period map $\kp^{\rm cpl}:\kt_\Gamma^{\rm cpl}\to Q_\Gamma$ is a covering and
thus, since $Q_\Gamma$ is simply connected, every connected component $\kt_0$ of $\kt_\Gamma^{\rm cpl}$ is generically mapped one-to-one onto $Q_\Gamma$.

\subsection{(Derived) Global Torelli Theorem}\label{derGT}
Before discussing the action of $\OO(\Gamma\oplus U)$ from the mirror symmetry
point of view we shall explain that a derived version of the Global Torelli
Theorem can be formulated by means of this action.

First, we reformulate the classical Global Torelli Theorem
for K3 surfaces as follows:

\begin{theorem}
Let $X$ and $X'$ be two K3 surfaces. Then $X$ and $X'$ are isomorphic if
and only if their images $\kp^{\rm cpl}(X,\varphi)$ and
$\kp^{\rm cpl}(X',\varphi')$ are contained in the same $\OO(\Gamma)$-orbit
in $Q_\Gamma$. 
\end{theorem}
(Of course, the choice of $\varphi$ and $\varphi'$ does not matter.)

In order to formulate the derived version of this, which consists in
weakening the isomorphism of $X$ and $X'$ to an equivalence of their
derived categories,  we need to complete the picture of the various
period maps as follows.

The diagram in Section \ref{defpmaps}
can be enriched by adding a moduli space that 
naturally contains the complex moduli space $\kt^{\rm cpl}_\Gamma$ and the
complex period domain $\Gr_2^{\rm po}(\Gamma_\IR)$ such that the group
$\OO(\Gamma\oplus U)$ acts naturally on the latter. 
We introduce the commutative diagram:
$$\xymatrix{\Gr_{2}^{\rm
po}(\Gamma_\IR)\times\Gamma_\IR\ar@{^{(}->}[r]^{\delta}&\Gr_{2}^{\rm
po}(\Gamma_\IR\oplus U_\IR)\\
\Gr_{2,1}^{\rm po}(\Gamma_\IR)\times\Gamma_\IR\ar@{>>}[u]
\ar@{^{(}->}[r]^{\gamma}&\Gr_{2,2}^{\rm
po}(\Gamma_\IR\oplus U_\IR)\ar@{>>}^{\pi}[u]}$$

Here, $\pi$ is the projection $(H_1,H_2)\mapsto H_1$ and $\delta$ is given
by $\delta:(P,B)\mapsto\{x-\langle x,B\rangle w~|~x\in P\}$.
This obviously yields the above commutative diagram. Moreover, $\pi$ is 
equivariant with respect to the natural $\OO(\Gamma\oplus U)$-action
on both spaces. But note that $\iota$ and $\tilde\xi=\iota\circ \xi$ do not
descend to $\Gr_2^{\rm po}(\Gamma_\IR\oplus U_\IR)$. 

Choosing a vanishing B-field for any marked K3 surface $(X,\varphi)$ yields
a map $\kt^{\rm cpl}_\Gamma\to\Gr_2^{\rm po}(\Gamma_\IR)\rpfeil{}{\delta}
\Gr_2^{\rm po}(\Gamma_\IR\oplus U_\IR)$.

Analogously to the discussion of the embedding $\gamma$ in Section 
\ref{PP}, one finds that the image of $\delta$ is not invariant
under the $\OO(\Gamma\oplus U)$-action, but it might of course happen that
the image of a marked K3 surface $(X,\varphi)$
 under some $\psi\in\OO(\Gamma\oplus U)
\setminus \OO(\Gamma)$ is mapped to the period of another K3 surface
$(X',\varphi')$. At least for algebraic K3 surfaces, when this happens can
be explained in terms of derived categories. This is due to a result
of Orlov \cite{Orlov} which is based on \cite{Mukai}.

\begin{theorem}
Two algebraic K3 surfaces $X$ and $X'$ have equivalent derived categories
$${\rm D}^{\rm b}({\rm Coh}(X))~ and~{\rm D}^{\rm b}({\rm Coh}(X'))$$ if and
only if their images $\delta\kp(X,\varphi)$ and $\delta\kp(X',\varphi')$ are
contained in the same $\OO(\Gamma\oplus U)$-orbit in
$\Gr_2^{\rm po}(\Gamma_\IR\oplus U_\IR)$.\qed
\end{theorem}

(Again, the choice of the markings $\varphi$ and $\varphi'$ is inessential.)

\begin{remark}
There is a conjecture that generalizes the above results to K3 surfaces with rational B-field $B\in H^2(X,\IQ)$.
The derived categories in this case have to be replaced by twisted derived categories, where one derives the abelian
category of coherent sheaves over an Azumaya algebra (cf.\ \cite{Caldararu}).
\end{remark}

The following result due to Hosono, Lian, Oguiso, Yau \cite{HLOY}
and independently to Ploog \cite{P} should be regarded as an analogue of the fact that the image of ${\rm Aut}(X)\to \OO^+(H^2(X,\IZ))$ is the set of Hodge
isometries. At the same time it is `mirror' to the result of Borcea discussed above.

\begin{theorem}
Let $X$ be a projective K3 surface. Then the image of $${\rm Autequ}(
{\rm D}^{\rm b}({\rm Coh}(X))\to\OO(H^*(X,\IZ))$$ contains the set of Hodge 
isometries contained in $\OO^+$.
\end{theorem}
Here the Hodge structure on $H^*(X,\IZ)$ is a weight-two Hodge structure
given by $H^{2,0}(X)\subset H^*(X,\IC)$.
As had been pointed out by Szendr${\rm\ddot{o}}$i in \cite{Szendroi}, mirror symmetry
suggests that the image should be exactly $\OO^{+}$. This would be the
analogue of Donaldson's result.


\section{Discrete group actions}\label{Discrete}

\medbreak All spaces considered in Section \ref{ps} are quotients
either of $\OO(\Gamma_\IR)$ or $\OO(\Gamma_\IR\oplus U_\IR)$. So
from a mathematical point of view it seems very natural to study
the action of the discrete groups $\OO(\Gamma)$ respectively
$\OO(\Gamma\oplus U)$ on these spaces. In fact, in order to obtain
moduli spaces of unmarked (complex) HK or (k{\"a}hler) IHS with or without B-fields.
one has to divide out by a smaller group. 
But in \cite{AM} it is argued that dividing out $\kt^{(4,4)}_\Gamma$
or $\kt_\Gamma^{(2,2)}$ by
$\OO(\Gamma\oplus U)$ yields the true moduli space of CFTs on K3
surfaces. In order to recover the full symmetry of the situation
they proceed as follows:

\medskip

{\bf i) Maximal discrete subgroups.}
Find a discrete group $G$ that acts on a certain moduli
space of relevant theories and show that it is maximal in the
sense that any bigger group would no longer act properly
discontinuously. (Recall that the quotient of a properly
discontinuous group action is  Hausdorff.)

{\bf ii) Geometric symmetries.}
Describe the part of $G$ (the geometric symmetries) that
identifies geometrically identical theories and the part that is
responsible for trivial identifications (e.g.\ integral shifts of
the B-field).

{\bf iii) Mirror symmetries.}
Show that $G$ is generated by the symmetries in {\bf
ii)} and a few others that are responsible for mirror symmetry
phenomena.

\subsection{Maximal discrete subgroups}
We first recall the following facts:

\medskip

$\bullet$ Let $G$ be a topological group which is Hausdorff and
locally compact. If $K\subset G$ is a compact subgroup then any
other subgroup $H$ acts properly discontinuously from the left on
the quotient space $G/K$ if and only if $H\subset G$ is a
discrete subgroup. (For the elementary proof see e.g.\ \cite[Lemma
3.1.1]{W}.)

$\bullet$ Let $L$ be a non-trivial definite even unimodular
lattice and let $q\geq3$. Then $\OO(L\oplus U^{\oplus q})\subset
\OO(L_\IR\oplus U_\IR^{\oplus q})$ is a maximal discrete subgroup
(cf.\ \cite{Allan}).

\bigskip

The second result in particular applies to the K3 surface lattice
$\Gamma=2(-E_8)\oplus 3U$ and yields that
$\OO(\Gamma)\subset\OO(\Gamma_\IR)$ and $\OO(\Gamma\oplus
U)\subset \OO(\Gamma_\IR\oplus U_\IR)$ are both maximal discrete
subgroups.

The group $\OO(\Gamma)$ acts on $\Gr_2^{\rm po}(\Gamma_\IR)$ and
$\Gr_3^{\rm po}(\Gamma_\IR)$. As we have seen
$$\Gr_2^{\rm
po}(\Gamma_\IR)\cong
\raisebox{.3ex}{$\OO(3,19)$}/\raisebox{-.6ex}{$\SO(2)\times\OO(1,19)$}~~{\rm
and}~~ \Gr_3^{\rm
po}(\Gamma_\IR)\cong\raisebox{.3ex}{$\OO(3,19)$}/\raisebox{-.6ex}{$\SO(3)\times\OO(19)$}.$$
In the second case we are in the above situation, i.e.\ the
quotient is taken with respect to the compact subgroup
$\SO(3)\times\OO(19)$. Hence, $\OO(\Gamma)$ acts properly
discontinuously on $\Gr_3^{\rm po}(\Gamma_\IR)$ and there is no
bigger subgroup of $\OO(\Gamma_\IR)$ than $\OO(\Gamma)$ with the
same property. However, the action of $\OO(\Gamma)$ on $\Gr_2^{\rm
po}(\Gamma_\IR)$ is badly behaved, as the subgroup
$\SO(2)\times\OO(1,19)$ is not compact. In fact, in the proof of Proposition
\ref{denseorb} we have already seen that the action of $\OO(\Gamma)$ is
not properly discontinuous.

We are more interested in the action of $\OO(\Gamma\oplus U)$ on
$\Gr_4^{\rm po}(\Gamma_\IR\oplus U_\IR)\cong
\OO(4,20)/(\SO(4)\times\OO(20))$. Again $\OO(\Gamma\oplus U)$ is
maximal discrete and $\SO(4)\times\OO(20)$ is compact.
Hence, there is no bigger properly
discontinuous subgroup action on $\Gr_4^{\rm po}(\Gamma_\IR\oplus
U_\IR)$. Analogously, one
finds that $\OO(\Gamma\oplus U)$ is a maximal discrete subgroup
of $\OO(\Gamma_\IR\oplus U_\IR)$ acting properly discontinuously
on $\Gr_{2,2}^{\rm po}(\Gamma_\IR\oplus U_\IR)$.

Presumably, all these arguments also apply to any HK manifold,
but details need to be checked. (Recall that $(H^2(X,\IZ),q_X)$)
is not necessarily unimodular in higher dimensions.)

\subsection{Geometric symmetries}\label{geomsym}

We will try to identify
``geometric'' symmetries and integral shifts of the B-field inside
$\OO(\Gamma\oplus U)$. To this end we use the identification
$$\phi:\Gr^{\rm po}_3(\Gamma_\IR)\times\IR_{>0}\times\Gamma_\IR\cong \Gr_4^{\rm po}(\Gamma_\IR\oplus U_\IR)$$
described in Section \ref{34spaces}.

The natural inclusion
$\OO(\Gamma)\subset\OO(\Gamma\oplus U)$ is compatible with this
is isomorphism, i.e.\ if $\Pi=\phi(F,\alpha,B)$ and
$\varphi\in\OO(\Gamma)\subset\OO(\Gamma\oplus U)$, then
$\varphi(\Pi)=\phi(\varphi(F),\alpha,\varphi(B))$. This is a
straightforward calculation which we leave to the reader.
Clearly, $\OO(\Gamma)$ acts naturally on all spaces $\kt$,
$\kt^{\rm met}$, $\kt^{\rm cpl}$, $\kt^{(2,2)}$, and
$\kt^{(4,4)}$ and the period maps are equivariant. Thus,
$\OO(\Gamma)$ is the subgroup that identifies geometrically
equivalent theories.

\medskip

Next let $B_0\in\Gamma$ and let
 $\varphi_{B_0} \in \OO(\Gamma \oplus U)$ be the automorphism
$w\mapsto w$, $w^\ast \mapsto B_0+w^\ast -
(B_0^2/2)w$, and $x\in \Gamma \mapsto x-\langle B_0,x\rangle w$. One easily
verifies that this really defines an isometry. We claim that
if $\Pi=\phi(F,\alpha,B)$, then $\varphi_{B_0}(\Pi)=\phi(F,\alpha, B+B_0)$.

In order to do this let us more generally consider an element
$\varphi\in\OO(\Gamma\oplus U)$ such that $\varphi(w)=w$. For
$\Pi\in \Gr_4^{\rm po}(\Gamma_\IR\oplus U_\IR)$, let
$\tilde\Pi:=\varphi(\Pi)$. Then $\tilde F'=\tilde\Pi\cap
w^\perp=\varphi(\Pi)\cap \varphi(w)^\perp=\varphi(\Pi\cap
w^\perp)=\varphi(F')$. Moreover, one has the two orthogonal
splittings $\tilde\Pi=\tilde B'\IR\oplus\tilde F'$ and $\tilde\Pi
=\varphi(B')\IR\oplus\varphi(F')$, where $\tilde B'$ is
determined by $\langle \tilde B',w\rangle=1$. Since
$\langle\varphi(B'),w\rangle=\langle\varphi(B'),\varphi(w)\rangle=\langle
B',w\rangle=1$, one concludes $\tilde B'=\varphi(B')$. In
particular, $\tilde B'^2=B'^2$. The B-field $B$ is given by
$B'=\alpha w+w^*+B$. Hence, $\tilde B'=\alpha
w+\varphi(w^*)+\varphi(B)$ and thus the B-field determined by
$\tilde B'$ is nothing but $\varphi(B)$.

All this applied to $\varphi=\varphi_{B_0}$ one finds that under the
isomorphism $\Gr_4^{\rm po}(\Gamma_\IR\oplus U_\IR)=\Gr_3^{\rm
po}(\Gamma_\IR)\times\IR_{>0}\times\Gamma_\IR$ the integral
B-shift by $B_0$ that maps $(F,\alpha,B)$ to $(F,\alpha, B+B_0)$
corresponds to $\varphi_{B_0}$.

We leave it to the reader to verify that also the
$\OO(\Gamma\oplus U)$-action on $\Gr^{\rm
po}_{2,2}(\Gamma_\IR\oplus U_\IR)$ is well-behaved in the sense
that $\OO(\Gamma)\subset\OO(\Gamma\oplus U)$ and the maps
$\varphi_{B_0}$ for $B_0\in\Gamma$ act on the subspace
$\gamma(\Gr_{2,1}^{\rm
po}(\Gamma_\IR)\times\Gamma_\IR)\subset\Gr_{2,2}^{\rm po}
(\Gamma_\IR\oplus U_\IR)$ in the natural way.



\subsection{Mirror symmetries}

The next result (due to C.\ T.\ C.\ Wall, \cite{Wall})
explains which additional group elements have to
be added in order to pass from $\OO(\Gamma)$ to $\OO(\Gamma\oplus
U)$.

\begin{proposition}\label{Wall}
Let $\Gamma$ be a unimodular lattice of index $(k,\ell)$ with
$k,\ell\geq2$. Then $\OO(\Gamma\oplus U)$ is generated by the
following three subgroups: $$~\OO(\Gamma),~ ~ \OO(U),~ {\rm and}
~~ \{\varphi_{B_0}~|~B_0\in\Gamma\}.$$\qed
\end{proposition}

Thus, the result applies to the K3 surface lattice $2(-E_8)\oplus
3U$, but presumably something similar can be said for the
case of the lattice $2(-E_8)\oplus 3U\oplus 2(n-1)\IZ$, which is realized by the
Hilbert scheme of a K3 surface.

In \cite{AM} passing from $\OO(\Gamma)$ to $\OO(\Gamma\oplus U)$
is motivated on the base of physical insight. As usual in
mathematical papers on mirror symmetry we will take this for
granted and rather study the effects of these additional
symmetries in geometrical terms. Thus, the rest of this paragraph
is devoted to the study a few special elements of $\OO(\Gamma\oplus U)$
that are not
contained in the subgroup generated by $\OO(\Gamma)$ and
$\{\varphi_{B_0}\ |\ B_0\in\Gamma\}$. In particular, we will be
interested in their induced action on $\Gr_{2,1}^{\rm
po}(\Gamma_\IR)\times\Gamma_\IR$.

So far we have argued that $\OO(\Gamma\oplus U)$ is a maximal
discrete subgroup of $\OO(\Gamma_\IR\oplus U_\IR)$ that acts on
the two period spaces that interest us: $\Gr_{2,2}^{\rm
po}(\Gamma_\IR\oplus U_\IR)$ and $\Gr_{4}^{\rm
po}(\Gamma_\IR\oplus U_\IR)$. However, there seems to be a bigger
group which naturally and properly discontinuously acts on the space $\Gr_{2,2}^{\rm
po}(\Gamma_\IR\oplus U_\IR)$ (which thus cannot be realized as a
subgroup of $\OO(\Gamma_\IR\oplus U_\IR)$).

\begin{definition}
The group $\tilde\OO(\Gamma\oplus U)$ is the group acting on
$\Gr_{2,2}^{\rm po}(\Gamma_\IR\oplus U_\IR)$ which is generated
by $\OO(\Gamma\oplus U)$ and the involution
$\iota:(H_1,H_2)\mapsto (H_2,H_1)$.
\end{definition}
Here $\bar H$ is the space $H$ with the opposite orientation.
Note that one could actually go further and consider the maps
$(H_1,H_2)\mapsto (H_1,\bar H_2)$ or $( H_1,H_2)\mapsto (\bar
H_1, H_2)$. However, for the versions of mirror symmetry that
will be discussed in these lectures $\iota$ will do.

\bigskip

Before turning to the mirror map that interests us most in Section
\ref{honestms} let us discuss a few more elementary cases:

\bigskip

\bigskip

\noindent{${{\bf -{\rm \bf id}_{U}}}$}\label{psi0}

\medskip

\noindent
Consider the automorphism $\psi_0\in\OO(\Gamma\oplus U)$ that
acts trivially on $\Gamma$ and as $-\id$ on $U$.

\begin{lemma}
The automorphism $\psi_0$ preserves the subspace $\Gr_{2,1}^{\rm
po}(\Gamma_\IR)\times\Gamma_\IR$ and acts on it by
$$((P,\omega),B)\mapsto((P,-\omega),-B).$$
\end{lemma}

\prf If $(H_1,H_2)=\gamma((P,\omega),B)$, then by definition
of $\psi_0$:

$$\psi_0(H_1)=\{x+\langle x,B\rangle w\ |\ x\in P\}=\{x-\langle x, (-B)\rangle w\ |\ x\in P\}$$
and
\begin{eqnarray*}
\psi_0(H_2)&=&(\frac{1}{2}
(\alpha-B^2)(-w)-w^*+B)\IR\oplus(\omega+\langle\omega,B\rangle w)\IR\\
&=&-(\frac{1}{2}(\alpha-(-B)^2)w+w^*-B)\IR\oplus(\omega-\langle\omega,(-B)\rangle w)\IR\\
\end{eqnarray*}
Thus, the sign of $\omega$ has to be changed in order to get the correct
orientation $\psi_0(H_2)$.\qed

\bigskip

\bigskip

\noindent{${\bf w\leftrightarrow w^*}$}\label{psi1}

\medskip

\noindent
Consider the automorphism $\psi_1\in\OO(\Gamma\oplus U)$ that
acts trivially on $\Gamma$ and by $\psi_1(w)=w^*$, $\psi_1(w^*)=w$
on $U$.

\begin{lemma}
The automorphism $\psi_1$ preserves the subspace
$\{((P,\omega),B)\ |\ B\in(P,\omega)^\perp,\alpha\ne B^2\}$ of
$\Gr_{2,1}^{\rm po}(\Gamma_\IR)\times\Gamma_\IR$ and acts on it by
$$((P,\omega),B)\mapsto\frac{2}{\alpha-B^2}((P,\omega),B).$$
\end{lemma}

\prf Indeed, by definition of $\psi_1$ one has $\psi_1(H_1)=H_1$ and
\begin{eqnarray*}
\psi_1(H_2)&=&\left(\frac{1}{2}(\alpha-B^2)w^*+w+B\right)
\IR\oplus(\omega-\langle\omega,B\rangle
w^*)\IR\\
&=&\left(w^*+\frac{2}{\alpha-B^2}w+\frac{2}{\alpha-B^2}B\right)\IR\oplus
\left(\frac{2}{\alpha-B^2}\omega\right)\IR\\
\end{eqnarray*}
Then check that for
$\tilde\omega:=\frac{2}{\alpha-B^2}\omega$
and $\tilde B:=\frac{2}{\alpha-B^2}B$ one indeed has
$\frac{2}{\alpha-B^2}=\frac{1}{2}(\tilde\omega^2-\tilde B^2)$.\qed

\bigskip

It is interesting to observe that on the yet smaller subspace $\{((P,\omega),0)\}$
the automorphism $\psi_1$ acts by $(P,\omega)\to\frac{2}{{\omega^2}}(P,\omega)$.
In the geometric context this will be interpreted as inversion of the volume or, in physical language, T-duality.

\begin{remark}
Nahm and Wendland argue that $w\leftrightarrow w^*$ occurs as an automorphism of the orbifold $(2,2)$-SCFT
associated to a Kummer surface. Thus, it has to be added as a global symmetry to the subgroup 
$\langle \OO(\Gamma),\{\varphi_B  ~| ~B\in\Gamma\}\rangle$. Due to the result of Wall, one thus obtains
the full $\OO(\Gamma\oplus U)$-action on $\Gr_4^{{\rm po}}(\Gamma_\IR\oplus U_\IR)$.

Note that in the original argument Aspinwall and Morrison had used another additional symmetry.
Writing $\Gamma\oplus U=(-E_8\oplus 2U)\oplus(-E_8\oplus 2U)$ allows one to consider the interchange
of the two summands $(-E_8\oplus 2U)\leftrightarrow (-E_8\oplus 2U)$ as an element in $\OO(\Gamma\oplus U)$.
This additional automorphism, which together with $\langle \OO(\Gamma),\{\varphi_B  ~| ~B\in\Gamma\}\rangle$ also
generates the whole $\OO(\Gamma\oplus U)$-action on $\Gr_4^{{\rm po}}(\Gamma_\IR\oplus U_\IR)$, is realized as an automorphism
of a certain Gepner model. For the details of both approaches we have to refer to the original articles.
\end{remark}
\subsection{The mirror map ${\bf U\leftrightarrow U'}$}\label{honestms}

If the lattice can be written as $\Gamma=\Gamma'\oplus U'$, where
$U'$ is a copy of the hyperbolic plane $U$, then by Wall's result
Proposition \ref{Wall}
the group $\tilde\OO(\Gamma\oplus U)$
is generated by $\OO(\Gamma)$, $\{\varphi_{B_0}\ |\
B_0\in\Gamma\}$, the involution $\iota$, and
$\xi\in\OO(\Gamma\oplus U)$ which is the identity on $\Gamma'$
and switches $U$ and $U'$. Here we use an isomorphism $U\cong U'$
which we fix once and for all. We consider $\Gr_{2,1}^{\rm po}(\Gamma_\IR)\times\Gamma_\IR$ as a  subspace of $\Gr_{2,2}^{\rm po}(\Gamma_\IR\oplus U_\IR)$ via
the injection $\gamma$.

Neither $\iota$ nor $\xi$ leave the subspace
$\Gr_{2,1}^{\rm po}(\Gamma_\IR)\times\Gamma_\IR$
invariant. Indeed, if $((P,\omega),B)$ then $H_1\subset
\Gamma_\IR\oplus\IR w$ and $H_2\not\subset \Gamma_\IR\oplus \IR
w$ and therefore $(H_2,H_1)=\iota(H_1,H_2)$ cannot be contained
in the image of $\gamma$. Similarly, for a general $(H_1,H_2)$
the pair of planes $(\xi(H_1),\xi(H_2))$ will not satisfy
$\xi(H_1)\subset\Gamma_\IR\oplus \IR w$.

\begin{definition}
$\tilde\xi:=\iota\circ\xi\in\tilde O(\Gamma\oplus U)$.
\end{definition}

By definition, $\tilde\xi$ acts naturally on $\Gr_{2,2}^{\rm
po}(\Gamma_\IR\oplus U_\IR)$ and $\Gr_4^{\rm po}(\Gamma_\IR\oplus
U_\IR)$. The action on the latter coincides with the action of
$\xi$.
 We will show that $\tilde\xi$ can
be used to identify certain subspaces of $\Gr_{2,1}^{\rm
po}(\Gamma_\IR)\times\Gamma_\IR$, but the whole $\Gr_{2,1}^{\rm
po}(\Gamma_\IR)\times\Gamma_\IR$ will again not be invariant.
Maybe it is worth emphasizing that $\tilde\xi$ is an involution.
Indeed, $\iota$ commutes with the action of $\OO(\Gamma_\IR\oplus
U_\IR)$ and both transformations
$\iota$ and $\xi$ are of order two.

Note that different decompositions  $\Gamma=\Gamma'\oplus U'$ yield different $\xi$, which then
relate different pairs of subspaces of
$\Gr_{2,1}^{\rm po}(\Gamma_\IR)\times\Gamma_\IR$. The following easy lemma
shows that we dispose of such a decomposition whenever we find
a hyperbolic plane contained in $\Gamma$.

\begin{lemma}\label{decomp}
If $U'$ is a hyperbolic plane contained in a lattice $\Gamma$, then $\Gamma=U'^\perp
\oplus U'$.
\end{lemma}

\prf Choose a basis $(v,v^*)$ of $U'$ that corresponds to the
basis $(w,w^*)$ of $U$ under the identification $U'\cong U$.
Furthermore, let $\Gamma':=U'^\perp$ and let $V$ be the subspace
of the $\IQ$-vector space $\Gamma_\IQ$ that is orthogonal to
$U'_\IQ$. Thus, $\Gamma_\IQ=V\oplus U'_\IQ$. Clearly,
$\Gamma'\subset V$ and, conversely, for any $v\in V$ there exists
$\lambda\in\IQ^*$ with $\lambda v\in V\cap\Gamma\subset\Gamma'$.
Hence, $V=\Gamma'_\IQ$. Let $x\in\Gamma$ and write $x=y+(\lambda
v+\mu v^*)$ with $y\in V$ and $\lambda,\mu\in\IQ$. Then $\langle
x,v\rangle ,\langle x,v^*\rangle\in\IZ$ implies
$\lambda,\mu\in\IZ$ and, therefore, $y=x-(\lambda v+\mu
v^*)\in\Gamma\cap V=\Gamma'$. Thus, $\Gamma=\Gamma'\oplus U'$.\qed

\bigskip

For the rest of this section we fix the orthogonal splitting
$\Gamma=\Gamma'\oplus U'$ together with an identification $U'=U$.
By ${\rm pr}:\Gamma_\IR\to\Gamma_\IR'$we denote the
orthogonal projection.

\begin{proposition}\label{msmap}
Let $((P,\omega),B)\in\Gr_{2,1}^{\rm
po}(\Gamma_\IR)\times\Gamma_\IR$ such that
$\omega,B\in\Gamma'_\IR\oplus\IR v$. Then the $\tilde\xi$-mirror
image $((P\dual,\omega\dual),B\dual):=\tilde\xi((P,\omega),B))$ is
again contained in $\Gr_{2,1}^{\rm
po}(\Gamma_\IR)\times\Gamma_\IR$. It is explicitly given as
\begin{eqnarray*}\sigma\dual&:=&\displaystyle{\frac{1}{\langle\re(\sigma),v\rangle}\left({\rm
pr}(B+i\omega)-\frac{1}{2}(B+i\omega)^2v+v^*\right)}\\
B\dual+i\omega\dual&:=&\displaystyle{\frac{1}{\langle\re(\sigma),v\rangle}\left({\rm
pr}(\sigma)-\langle\sigma,B\rangle v\right)}\\
\end{eqnarray*}
Here, we have replaced $P$ by the corresponding line $[\sigma]\in
Q_\Gamma\subset\IP(\Gamma_\IC)$. Furthermore, we have chosen
$\sigma$  such that $\im(\sigma)$ is orthogonal to $v$.
\end{proposition}

\prf By definition the positive plane $P$ is contained in
$\omega^\perp$. Since the intersection of $\omega^\perp$ with
$\Gamma'_\IR\oplus \IR v$ and $\Gamma'_\IR\oplus\IR v^*$ have
both only one positive direction, $P$ cannot be contained in
either of them. Thus, we may choose $\sigma$ such that
$v^\perp\cap P=\im(\sigma)\IR$ and
$\langle\re(\sigma),v\rangle\ne0$. This justifies the above
choices. Also note that $\omega\dual$ and $B\dual$ do not change
when $\sigma$ is changed by a real scalar. The defining equations
for $B\dual+i\omega\dual$ and $\sigma\dual$ are spelled out as
follows
\begin{eqnarray*}
\sigma\dual&:=&\displaystyle{\frac{1}{\langle\re(\sigma),v\rangle}\left(-\frac{1}{2}(B+i\omega+v^*)^2v+B+i\omega+v^*\right)}\\
\omega\dual&:=&\displaystyle{\frac{1}{\langle\re(\sigma),v\rangle}\left(\im(\sigma)-\langle\im(\sigma),v^*\rangle v-\langle\im(\sigma),B\rangle v\right)}\\
B\dual&:=&\displaystyle{\frac{1}{\langle\re(\sigma),v\rangle}\left(\re(\sigma)-\langle\re(\sigma),v\rangle
v^*-\langle\re(\sigma),v^*\rangle v-\langle\re(\sigma),B\rangle v\right)}\\
\end{eqnarray*}
Let us now compute $\tilde\xi(H_1,H_2)$. We denote
$\gamma((\sigma\dual,\omega\dual),B\dual)$ by
$(H_1\dual,H_2\dual)$, where $\sigma\dual,\omega\dual$, and
$B\dual$ are as above.

The space $H_1\dual$ is spanned by the real and imaginary part of
$\sigma\dual-\langle\sigma\dual,B\dual\rangle w$. A simple
calculation yields
$$\langle\sigma\dual,B\dual\rangle=-\langle\re(\sigma),v\rangle^{-1}\left(\langle
B,v^*\rangle+i\langle\omega,v^*\rangle\right).$$ Thus, $H_1\dual$
is spanned by
\begin{eqnarray*}
&&B+v^*+\frac{1}{2}(\omega^2-B^2-2\langle B,v^*\rangle) v+\langle
B,v^*\rangle w\\
&=&\frac{1}{2}(\omega^2-B^2)v+v^*+(B-\langle B,v^*\rangle
v)+\langle B,v^*\rangle w\\
&=&\xi\left(\frac{1}{2}(\omega^2-B^2)w+w^*+(B-\langle B,
v^*\rangle v)+\langle B,v^*\rangle v\right)\\
&=&\xi\left(\frac{1}{2}(\omega^2-B^2)w+w^*+B\right)
\end{eqnarray*}
and
\begin{eqnarray*}
&&\omega-\langle \omega,B+v^*\rangle v+\langle \omega,v^*\rangle w\\
&=&(\omega-\langle\omega,v^*\rangle v)+\langle\omega,v^*\rangle
w-\langle\omega,B\rangle v\\
&=&\xi\left(\omega-\langle\omega,v^*\rangle
v+\langle\omega,v^*\rangle v-\langle\omega,B\rangle w\right)\\
&=&\xi\left(\omega-\langle\omega,B\rangle w\right).
\end{eqnarray*}
Hence, $H_1\dual=\xi(H_2)$. Similarly, one proves
$H_2\dual=\xi(H_1)$. First one computes
\begin{eqnarray*}
{\omega\dual}^2&=&\langle\re(\sigma),v\rangle^{-2}\im(\sigma)^2,\\
{B\dual}^2&=&\langle\re(\sigma),v\rangle^{-2}\left(\re(\sigma)^2-2\langle\re(\sigma),v\rangle\langle\re(\sigma),v^*\rangle\right)\\
\langle\omega\dual,B\dual\rangle&=&-\langle\re(\sigma),v\rangle^{-1}\langle\im(\sigma),v^*\rangle,
\end{eqnarray*}
where one uses $\langle\im(\sigma),v\rangle=0$. Since
$\im(\sigma)^2=\re(\sigma)^2$, this yields
$${\omega\dual}^2-{B\dual}^2=2\langle\re(\sigma),v\rangle^{-1}\langle\re(\sigma),v^*\rangle.$$
Hence, $H_2\dual$ is spanned by
\begin{eqnarray*}
&&\frac{1}{2}({\omega\dual}^2-{B\dual}^2)w+w^*+B\dual\\
&=&\langle\re(\sigma),v\rangle^{-1}\langle\re(\sigma),v^*\rangle
w+w^*\\
&&+\langle\re(\sigma),v\rangle^{-1}\left(\re(\sigma)-\langle\re(\sigma),v\rangle
v^*\langle\re(\sigma),v^*\rangle v-\langle\re(\sigma),B\rangle
v\right)
\end{eqnarray*}
and
\begin{eqnarray*}
&&{\omega\dual}-\langle\omega\dual,B\dual\rangle w\\
&=&\langle\re(\sigma),v\rangle^{-1}\left((\im(\sigma)-\langle\im(\sigma),v^*\rangle
v)-\langle\im(\sigma),B\rangle v+\langle\im(\sigma),v^*\rangle
w\right)
\end{eqnarray*}
Thus, $\xi(H_2\dual)$ is generated by
$\re(\sigma)-\langle\re(\sigma),B\rangle w$ and
$\im(\sigma)-\langle\im(\sigma),B\rangle w$. Hence,
$H_2\dual=\xi(H_1)$.\qed

\bigskip

\begin{examples} The proposition can be used to identify certain
subspaces of $\Gr_{2,1}^{\rm po}(\Gamma_\IR)\times\Gamma_\IR$ via
the mirror map $\tilde\xi$. We will present a few
examples, which will be interpreted geometrically later on. As the B-field
from a geometric point of view is not well-understood,
 we will be especially interested in those points with vanishing B-field.

{\bf i)} Fix an orthogonal decomposition $\Gamma'_\IR =V\oplus
V\dual$, such that both subspaces $V$ and $V\dual$ contain a
positive line. The automorphism $\tilde\xi\in\tilde
O(\Gamma\oplus U)$ induces a bijection between the two subspaces:
$$\{((P,\omega),B)\ |\ B,\omega\in V,\ P\subset V\dual\oplus U'_\IR\}~{\rm
and}~\{((P,\omega),B)\ |\ B,\omega\in V\dual,\ P\subset V\oplus
U'_\IR\}.$$ Note that in this case the formulae for
$(\sigma\dual,\omega\dual,B\dual)$ simplify slightly to:
$$\sigma\dual=\frac{1}{\langle\re(\sigma),v\rangle}(-\frac{1}{2}(B+i\omega)^2v+v^*+B+i\omega),~
\omega\dual=\frac{1}{\langle\re(\sigma),v\rangle}(\im(\sigma)-\langle\im(\sigma),v^*\rangle
v),$$ and
$$B\dual=\frac{1}{\langle\re(\sigma),v\rangle}(\re(\sigma)-\langle\re(\sigma),v\rangle
v^*-\langle\re(\sigma),v^*\rangle v).$$
From here it is easy to
verify that $\tilde\xi$ maps these two subspaces into each
other. Note that $B\dual+i\omega\dual$ is  up to the scalar
factor $\langle\re(\sigma),v\rangle^{-1}$ nothing but the
projection of $\sigma\in V\dual_\IC\oplus U'_\IC$ to $V\dual_\IC$.
\medskip

{\bf ii)} It might be interesting to see what happens in the
previous example if we set the B-field zero. Under the assumption
of {\bf i)} the symmetry $\tilde\xi$ induces a bijection
between the following two subspaces
\begin{eqnarray*}
&\{((P,\omega),B=0)~|~\omega\in V, ~\re(\sigma)\in U'_\IR,
~\im(\sigma)\in V\dual\}\\
{\rm and}&~\{((P,\omega),B=0)~|~\omega\in V\dual, ~\re(\sigma)\in
U'_\IR, ~\im(\sigma)\in V\}
\end{eqnarray*}
Indeed, $\im(\sigma\dual)=
\langle\re(\sigma),v\rangle^{-1}(-\langle B,\omega\rangle
v+\omega)$ and
$\re(\sigma\dual)=\langle\re(\sigma),v\rangle^{-1}(\frac{1}{2}
(\omega^2-B^2)v+v^*+B)$. Thus, if $B=0$ one has
$\im(\sigma\dual)=\langle\re(\sigma), v\rangle^{-1}\omega$ and
$\re(\sigma\dual)=\langle\re(\sigma),
v\rangle^{-1}(\frac{\omega^2}{2}v+v^*)$, and $B\dual=0$.
Conversely, if $\re(\sigma\dual)\in U'_\IR$ then $B=0$. Moreover,
$\im(\sigma)\in V$ implies $\omega\in V$ and $\omega\dual\in
V\dual$ implies $\im(\sigma)\in V\dual$. Eventually, $B\dual=0$
yields $\re(\sigma)\in U'_\IR$.

\medskip

{\bf iii)} In this example we will not need any further
decomposition of $\Gamma'_\IR$. The automorphism
$\tilde\xi\in\tilde \OO(\Gamma\oplus U)$ induces an involution
on the subspace
$$\{((P,\omega),B)~|~\omega,B\in\Gamma_\IR'\oplus\IR
v\}\subset\Gr_{2,1}^{\rm po}(\Gamma_\IR).$$ This follows again
easily from the explicit description of
$(\sigma\dual,\omega\dual,B\dual)$.

\medskip

{\bf iv)} Also in {\bf iii)} one finds a smaller subset
parametrizing only objects with trivial B-field that is left invariant
by $\tilde\xi$. Indeed, the subspace
$$\{((P,\omega),0)~|~\omega\in\Gamma_\IR',~P\cap U'_\IR\ne0\}$$
is mapped onto itself under $\tilde\xi$.\qed
\end{examples}

\begin{remark}  If $((P,\omega),B)$ such that $B\dual=0$ and
$\varphi\in\OO(\Gamma')$, then also
$\tilde\xi((\varphi(P),\varphi(\omega)),\varphi(B))$ has
vanishing B-field. Geometrically this is used to argue that if
the mirror $X\dual$ of $X$ has vanishing B-field then the same
holds for the mirror of $f^*X$ under any diffeomorphism $f$ of
$X$ with $f^*|_{U'}={\rm id}$.
The assertion is an immediate consequence of the explicit description of
$B\dual$ given above (cf.\ \cite{Szendroi}).
\end{remark}

 Note that $\tilde\xi$ is by far the most interesting
automorphism considered so far, as it really mixes the `complex
direction' $\sigma$ with the `metric direction' $(\omega,B)$.
However, at least for the case of the K3 lattice
$\Gamma=2(-E_8)\oplus 3U$ the automorphisms $\xi$ respectively
$\{-\id_U,w\leftrightarrow w^*\}$ together with
$\{\OO(\Gamma),\varphi_{B_0\in\Gamma}\}$ generate both the same
group, namely $\OO(\Gamma\oplus U)$. This is a consequence of
Proposition \ref{Wall}, where one uses $\xi\OO(U)\xi=\OO(U')$ and thus
$\OO(U)\subset\langle\xi,\OO(\Gamma)\rangle$. So in this sense,
$\xi\in\OO(\Gamma\oplus U)$ as an automorphism of $\Gr_{2,2}^{\rm
po}(\Gamma_\IR\oplus U_\IR)$ is not more or less interesting than
those in \ref{psi0} and \ref{psi1}, but for the latter ones the
interesting things happen outside the `geometric world' of
$\Gr_{2,1}^{\rm po} (\Gamma_\IR)\times\Gamma_\IR$.


\section{Geometric interpretation of mirror symmetry}\label{msK3}

\subsection{Lattice polarized mirror symmetry}

Let $\Gamma$ as before be the K3 lattice $2(-E_8)\oplus 3U$ and
fix a sublattice $N\subset \Gamma$ of signature $(1,r)$.

\begin{definition}
An $N$-polarized marked K3 surface is a marked K3 surface $(X,\varphi)$
such that $N\subset\varphi({\rm Pic}(X))$.
\end{definition}

Note that any $N$-polarized K3 surface is projective. If $\kt^{\rm
cpl}_\Gamma$ is the moduli space of marked K3 surfaces we denote
by $\kt_{N\subset\Gamma}^{\rm cpl}$ the subspace that consists of
$N$-polarized marked K3 surfaces. Analogously, one defines
$$\kt^{(2,2)}_{N\subset\Gamma}\subset\kt^{(2,2)}_\Gamma$$
as the subset of all marked K{\"a}hler K3 surfaces with B-field
$(X,\omega,B,\varphi)$ such that $N\subset{\rm Pic}(X)$ and
$\omega,B\in N_\IR$. Here and in the following, we omit the
marking in the notation, i.e.\ the identification $H^2(X,\IZ)\cong
\Gamma$ via $\varphi$ will be understood.

The condition $N\subset{\rm Pic}(X)$ is in fact equivalent to
$V:=N_\IR\subset{\rm Pic}(X)_\IR$. The latter can furthermore be rephrased
as $V\subset (H^{2,0}(X)\oplus H^{0,2}(X))^\perp$, i.e.\
$\sigma\in V^\perp_\IC$.

By construction there exists a natural map
$$\kt_{N\subset\Gamma}^{(2,2)}\to\kt_{N\subset\Gamma}^{\rm cpl}.$$
The fibre over $(X,\varphi)\in\kt_{N\subset\Gamma}^{\rm cpl}$ is
isomorphic to $V_\IR+i({\cal K}_X\cap V_\IR)$ via
$(\omega,B)\mapsto B+i\omega$.

Using the period map, the space $\kt_{N\subset\Gamma}^{(2,2)}$ can
be realized as a subspace of $\Gr_{2,1}^{\rm
po}(\Gamma_\IR)\times\Gamma_\IR\subset\Gr_{2,2}^{\rm
po}(\Gamma_\IR\oplus U_\IR)$. Its closure
$\overline\kt_{N\subset\Gamma}^{(2,2)}$ consists of all points
 $((P,\omega),B)\in\Gr_{2,1}^{\rm
po}(\Gamma_\IR)\times\Gamma_\IR$ such that $P\subset V^\perp$ and
$\omega,B\in V$.
Indeed, via the period map $\kt^{(2,2)}_{N\subset\Gamma}$ is identified
with an open subset of $\{((P,\omega),B)~|~B,\omega\in V, P\subset
V^\perp\}$ and the latter is irreducible.

Let us now assume that the orthogonal complement
$N^\perp\subset\Gamma$ contains a hyperbolic plane $U'\subset
N^\perp$. Then $N^\perp=N\dual\oplus U'$ by Lemma \ref{decomp} for some
sublattice $N\dual\subset\Gamma$ of signature $(1,18-r)$. The
real vector space $N\dual_\IR$ is denoted by $V\dual$. As above
one introduces $\kt_{N\dual\subset\Gamma}^{(2,2)}$ and
$\kt_{N\dual\subset\Gamma}^{\rm cpl}$.

\begin{proposition}
The mirror symmetry map $\tilde\xi$ associated to the splitting
$\Gamma=\Gamma'\oplus U'$ induces a bijection
$$\overline\kt_{N\subset\Gamma}^{(2,2)}\cong\overline\kt_{N\dual\subset\Gamma}^{(2,2)}.$$
\end{proposition}

\prf
By the description of $\overline\kt^{(2,2)}_{N\subset\Gamma}$ as the set
$\{((P,\omega),B)~|~B,\omega\in V, P\subset V^\perp\}$, it suffices to show
that the mirror map identifies the two sets
$\{((P,\omega),B)~|~B,\omega\in V, P\subset V^\perp\}$ and
$\{((P,\omega),B)~|~B,\omega\in V\dual, P\subset (V\dual)^\perp\}$,
which has been observed already in the Examples in Section \ref{honestms}.
\qed

\begin{remark} {\bf i)} In general, we cannot expect to have
a bijection
$\kt_{N\subset\Gamma}^{(2,2)}\cong\kt_{N\dual\subset\Gamma}^{(2,2)}$.
Indeed, for a point in $\kt_{N\subset\Gamma}^{(2,2)}$ that
corresponds to a triple $((P,\omega),B)$ the image
$((P\dual,\omega\dual),B\dual)=\tilde\xi((P,\omega),B)$ might
admit a $(-2)$-class $c\in(P\dual)^\perp\cap\Gamma$ with $\langle
c,\omega\dual\rangle=0$. In fact, these two conditions on the
$(-2)$-class $c$ translate into the equations $\langle c+\langle
c,v\rangle B,\omega\rangle=0$ and $\langle c-\langle c,v\rangle
v^*,\im(\sigma)\rangle=0$. To exclude this possibility one would
need to derive from this fact that there exists a $(-2)$-class
$c'$ with $\langle c',\sigma\rangle=0$ and $\langle
c',\omega\rangle=0$ and this doesn't seem possible in general.

One should regard this phenomenon as a very fortunate fact. As
points in the boundary are interpreted as singular K3 surfaces,
it enables us to compare smooth K3 surfaces with singular ones.
One should try to construct examples of (singular) Kummer
surfaces in this context.

{\bf  ii)} Also note that if $\omega\in\Gamma$, i.e.\ $\omega$ corresponds
to a line bundle, then $\omega\dual$ does not necessarily have the
same property.

{\bf  iii)} We also remark that the lattices $N$ and $N\dual$ are
rather unimportant in all this. Indeed, what really matters are
the two decompositions $\Gamma=\Gamma\oplus U'$ and
$\Gamma'_\IR=V\oplus V\dual$.

\end{remark}

To conclude this section, we shall compare the above discussion with 
\cite{Dolgachev}. Let $N\subset\Gamma$ and $N^\perp=N\dual\oplus U'$ be as
before. Following \cite{Dolgachev} one defines
$$\Omega:=N\dual_\IR\oplus i(N\dual_\IR\cap\kc)~{\rm and}~
D_N:=Q_\Gamma\cap\IP(N_\IC^\perp).$$

Then by \cite[Thm.4.2,Rem.4.5]{Dolgachev} the map
$$\alpha:\Omega\to D_N,z\mapsto[z-\frac{1}{2}z^2\cdot v+v^*]$$
is an isomorphism. This map obviously coincides with
$(B+i\omega)\mapsto[\sigma\dual]$ as described in Proposition \ref{msmap}, since for
$B,\omega\in N_\IR\dual\subset\Gamma_\IR'$ one has ${\rm
  pr}(B+i\omega)=B+i\omega$. Thus, the map $\alpha$ coincides with the map
given by the isomorphism
$\overline\kt^{(2,2)}_{N\dual\subset\Gamma}\cong\overline\kt^{(2,2)}_{N\subset\Gamma}$.
To make this precise note that $\kt^{\rm
  cpl}_{N\subset\Gamma}\cong D_N$ via the period map and that
$((P,\omega),B)\mapsto B+i\omega$ defines a surjection
$\overline\kt_{N\dual\subset\Gamma}^{(2,2)}\epimorph\Omega$. This yields a
commutative diagram

$$\xymatrix{
\overline\kt^{(2,2)}_{N\dual\subset\Gamma}\ar @{>}[r]^{\cong}\ar@{>>}[d]&\overline\kt^{(2,2)}_{N\subset\Gamma}\ar@{>>}[d]\\
\Omega\ar@{>}[r]^{\cong~~~~~} &D_N\cong\kt^{\rm cpl}_{N\subset\Gamma}\\}
$$
which emphasizes the fact that the mirror isomorphism identifies K{\"a}hler
deformations with complex deformations.

\begin{remark}\label{LPrem}
{\bf i)} We also mention the following result of Looijenga and Peters 
\cite{LP}, which shows that lattices of small rank can always be realized.
Let $\Gamma$ be the K3 lattice and $N$ any even lattice of rank at most three.
Then there exists a primitive embedding $N\subset \Gamma$. If the rank is smaller
than three then this primitive embedding is unique up to automorphisms
of $\Gamma$, i.e.\ elements of $\OO(\Gamma)$. An even more general
version of this result can be found in \cite{Nikulin}.

{\bf ii)} Moduli spaces of polarized K3 surfaces and their compactifications
have been treated in detail by Looijenga, Friedman, Scatonne and many others
(see \cite{Scat}).
\end{remark}

\subsection{Mirror symmetry by hyperk{\"a}hler rotation}\label{MSbyHK}

Let $X$ be a K3 surface with K{\"a}hler class $\omega_I$. Assume that
a decomposition $H^2(X,\IZ)=\Gamma=\Gamma'\oplus U'$ together with
an isomorphism $\xi:U'\cong U$ has been fixed. As before, we denote by
$(v,v^*)$ the basis of $U'$ that corresponds to $(w,w^*)$.

\begin{proposition}\label{MSHK}
  Let $\omega_I\in\Gamma'_\IR$ be a K{\"a}hler class on $X$ and assume
  that $\sigma_I\in H^{2,0}(X)$ with $\sigma_I\bar\sigma_I=2\omega_I^2$
  can be chosen such that ${\rm Re}(\sigma_I)\in U'_\IR$,
  $\im(\sigma_I)\in\Gamma'_\IR$, and
  $\langle\re(\sigma_I),v\rangle=1$. Then the
$\tilde\xi$-mirror of $(X,\omega_I,B=0)$ is given by the formula
$$\begin{array}{rcl}
\sigma\dual&:=&\displaystyle{\frac{1}{\langle\re(\sigma_I),v\rangle}\left(\re(\sigma_I)+i\omega_I\right)}\\
\omega\dual&:=&\displaystyle{\frac{1}{\langle\re(\sigma_I),v\rangle}\im(\sigma_I)~{\rm and}~ B\dual=0}\\
\end{array}$$
\end{proposition}

\prf Using $\im(\sigma_I)\in\Gamma_\IR'$ and the general formula given in the proof
of Proposition \ref{msmap}, we find that it is enough to prove $\re(\sigma_I)=(\omega_I^2/2)v+v^*$, but this follows immediately from the assumption
$\re(\sigma_I)\in U_\IR'$ and $\re(\sigma_I)^2=\omega_I^2$.\qed

\bigskip

Note that, if we only know that $\re(\sigma_I)\in U_\IR$, the
metric, and hence $\sigma_I$ and $\omega_I$, can be rescaled such that
$\langle(\sigma_I),v\rangle=1$. But scaling the metric changes the
complex structure of the mirror. This will be important when we
discuss the large K{\"a}hler and complex structure limits.
\begin{corollary}
If $(X,\omega_I)$ is a K3 surface as in the proposition then the mirror K3 surface $X\dual$ is
obtained by hyperk{\"a}hler rotation to $-K$.
\end{corollary}

\prf Observe that $\re(\sigma_I)=\omega_J$. Thus, $\sigma\dual=
\langle\re(\sigma_I),v\rangle^{-1}(\omega_J+i\omega_I)$.
On the other hand, $\sigma_{-K}=\omega_J+i\omega_I$. Hence, 
$X\dual$ is given by the complex structure $-K$.\qed

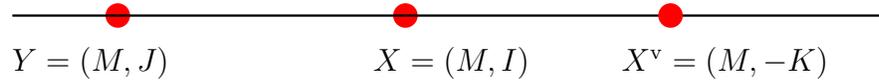
\psset{unit=1.5cm}
\hspace{2cm}
\begin{pspicture*}(0,1)(7.733333,3.583333)
\psset{dotstyle=x, dotscale=1.0, linewidth=0.02}
\psdots[dotscale=2.5,dotstyle=*,linecolor=red](0.933334,2.483334)
\put(0.0,2.0){$Y=(M,J)$}
\psdots[dotscale=2.5,dotstyle=*,linecolor=red](3.483333,2.483333)
\put(3.2,2.0){$X=(M,I)$}
\psline[linecolor=black](0.000000,2.483333)(7.733333,2.483333)
\psdots[dotscale=2.5,dotstyle=*,linecolor=red](5.833314,2.483333)
\put(5.4,2.0){$X\dual=(M,-K)$}
\end{pspicture*}

There is one tiny subtlety. If we compute also the mirror K{\"a}hler form
$\omega\dual$, we obtain $-\omega_{-K}$. But this is of no importance,
as we can always apply the harmless global transformation
$-{\rm id}\in\OO(\Gamma)$.

\begin{remark}\label{doubtsHKMS}
Thus, in the  very special case that $H^2(X,\IZ)=\Gamma'\oplus U'$ such that
$\re(\sigma_I)\in U_\IR'$ and $\omega_I\in\Gamma'_\IR$, mirror symmetry
is given by hyperk{\"a}hler rotation. However, the phenomenon seems
rather accidental and one should maybe not expect that there
is is a deeper interplay between mirror symmetry and hyperk{\"a}hler
rotation. E.g.\ one can check that the solution of
the Maurer-Cartan equation given by the Tian-Todorov coordinates
is not the one obtained from hyperk{\"a}hler rotation and deforming the
hyperk{\"a}hler structure.

\end{remark}
\subsection{Mirror symmetry for elliptic K3 surfaces}

Here we will discuss one geometric instance where special K3 surfaces as
treated in the last section naturally occur.

Let $\pi:Y\to \IP^1$ be an elliptic K3 surface with a section
$\sigma_0\subset Y$. The cohomology class $f$ of the fibre and $[\sigma_0]$ generate a sublattice $U'\subset H^2(Y,\IZ)$. It can be identified with
the standard hyperbolic plane by choosing as a basis
$v=f$ and $v^*=f+\sigma_0$. Thus, we obtain a decomposition
$\Gamma:=H^2(Y,\IZ)=\Gamma'\oplus U'$ together with an isomorphism $\xi:U'\cong U$,.

 Let us now study the action of $\tilde\xi$ on K3
surfaces that are related to $Y$. If we fix a HK-metric $g$ on
$Y$, then we may write $Y=(M,J)$, where $J$ is one of the
compatible complex structures $\{aI+bJ+cK~|~a^2+b^2+c^2=1\}$
associated with $g$. A holomorphic two-form on $Y$ can be given as
$\sigma_J=\omega_K+i\omega_I$. The reason why the complex
structure that defines $Y$ is denoted  $J$ is that we will
actually not describe the mirror of $Y$, but rather of $X:=(M,I)$.

Clearly $X$ inherits the torus fibration from $Y$ which gives rise
to a differentiable map $\pi:X\to \IP^1$.

\begin{lemma}
The torus fibration $\pi:X\to \IP^1$ is a SLAG fibration.
\end{lemma}

\prf Indeed, since the holomorphic two-form $\sigma_J$ vanishes
on any holomorphic curve in $Y$, the form $\omega_I=\im(\sigma_J)$ vanishes
in particular on every fibre of $X\to\IP^1$, i.e.\ all fibres are Lagrangian.
Moreover, since $\sigma_I=\omega_J+i\omega_K$ and $\omega_K=\frac{1}{2}(\sigma_J+\bar\sigma_J)$, we see that $\im(\sigma_I)|_{\pi^{-1}(t)}=0$ and
$\re(\sigma_I)|_{\pi^{-1}(t)}=\omega_J|_{\pi^{-1}(t)}$.
Hence, the (smooth) fibres are special Lagrangians of phase $0$.\qed

\bigskip

Again as a consequence of the general formula in Proposition \ref{msmap}
one computes the mirror of $(X,\omega_I)$ explicitly.

\begin{proposition}
The $\tilde\xi$-mirror of $(X,\omega_I)$ is the K3 surface
$X\dual$ given by the period
$$\sigma\dual=\frac{1}{{\rm vol}(f)}\left(\frac{\omega_I^2}{2}f+\sigma_0+f+i\omega_I\right),$$ which is endowed with the K{\"a}hler class
$$\omega\dual=\frac{1}{{\rm vol}(f)}\im(\sigma_I),$$
where ${\rm vol}(f)=\langle \omega_J,f\rangle$ is the volume of the fibre of the
elliptic fibration $Y\to\IP^1$.\qed
\end{proposition}

Of course, an explicit formula could also be given for the mirror of $X$
endowed with the K{\"a}hler form $\omega_I$ and an auxiliary B-field. We leave this
to the reader.

\begin{remark}
A priori, the mirror K3 surface could be singular,
i.e.\ there could be a $(-2)$-class $c\in \Gamma$ such that $\langle c,\omega\dual\rangle=\langle c,\sigma\dual\rangle=0$. For such a class we would have 
$\langle c,\omega_I\rangle=\langle c,\im(\sigma_I)\rangle=0$. Of course,
if we also had $\langle c,\omega_I\rangle=0$, then already $(X,\omega_I)$ would
be singular.
\end{remark}

If in addition we choose $\omega_J:=(\alpha/2+1)f+\sigma_0$ is a
K{\"a}hler class for some $\alpha>0$, e.g.\ when 
$Y$ has Picard number two, then $\re(\sigma_I)$ satisfies the condition of 
Proposition \ref{MSHK}, i.e.\ $\re(\sigma)\in U'_\IR$
$\im(\sigma_I)\in\Gamma'_\IR$, and $\langle\re(\sigma),v\rangle=1$. Hence, in this case
the mirror of $X$ is given by the complex structure $-K$
and the fibration is still a SLAG fibration.

\begin{remark}
If we go back to the more general case, where $\omega_J$ on $Y$ might
 be arbitrary, then we still see that $\omega\dual$, $\im(\sigma\dual)
\in\langle v,v^*\rangle^\perp$, i.e.\ at least cohomologically the classes
$f$ and $\sigma_0$ are still Lagrangian on the mirror, as in the more special
case above where $X\dual$ was given by $-K$.
\end{remark}


\subsection{FM transformation and mirror symmetry}

Here we shall explain how the mirror symmetry map $\tilde\xi$ can be viewed
as the action of the FM-transformation on cohomology. 
Let us first recall the setting of the previous section. Let $\pi:Y\to\IP^1$
be an elliptic K3 surface with a section $\sigma_0$. The cohomology class
of the fibre will be denoted by $f$ and the complex structure by $J$,
i.e.\ $Y=(M,J)$. Assume that $[\omega_J]=(\alpha/2+1)f+\sigma_0$
is a K{\"a}hler class on $Y$. Then let
$X=(M,I)$, where $I$ and $K=IJ$ are the other two compatible complex
structures associated to the hyperk{\"a}hler metric underlying $[\omega_J]$.
Then $\pi:X\to\IP^1$ is a SLAG fibration. Furthermore, consider the dual elliptic fibration $\pi:Y\dual\to\IP^1$. Since $\pi:Y\to\IP^1$ has a section, there
is a canonical isomorphism $Y\cong Y\dual$ compatible with the projection.

Let $\kp\to Y\times_{\IP^1}Y\dual=Y\times_{\IP^1} Y$ be the relative Poincar{\'e}
sheaf and let $\tilde{\rm FM}:H^*(Y,\IZ)\cong H^*(Y,\IZ)$ denote the
cohomological Fourier-Mukai transformation $\beta\mapsto
q_*(p^*\beta.{\rm ch}(\kp))$. A standard calculation shows (cf.\ \cite{Bruzzo}):

\begin{lemma}
$\tilde{\rm FM}=\xi$
\end{lemma}
\prf For the Fourier-Mukai transform on the level
of derived categories one has ${\rm FM}({\cal O}_f)={\cal O}_{\sigma_0\cap f}$ and
${\rm FM}({\cal O}_{\sigma_0}(-1))={\cal O}_Y$. Since for the Mukai vectors one has 
$v({\cal O}_f)=f$, $v({\cal O}_{\sigma_0}(-1))=\sigma_0$, $v({\cal O}_{\sigma_0\cap f})=
[{\rm pt}]\in H^4(Y)$, and $v({\cal O}_Y)=[{\rm pt}]+[Y]\in H^4(Y)\oplus H^0(Y)$,
passing to cohomology yields the result.
\qed

\bigskip

Thus, on a purely cohomological level, this fits nicely with the expectation
that SLAGs on $X$ should correspond to holomorphic objects on the
mirror $X\dual$ which happens to be $(M,-K)$ as was explained before.
Indeed, the two SLAGs $f$ and $\sigma_0$ on $X$ are first hyperk{\"a}hler rotated
to holomorphic objects in $Y$, namely the holomorphic fibre respectively section of the elliptic fibration $Y\to \IP^1$. The FM-transforms of those
are $k(\sigma_0\cap\pi^{-1}(t))$ respectively ${\cal O}_Y$.
Although we still have to hyperk{\"a}hler rotate from $Y\dual=Y$
to $X\dual=(M,-K)$ the
cycles $k(\sigma_0\cap\pi^{-1}(t))$ and ${\cal O}_Y$ stay holomorphic. Thus,
the mirror symmetry map $\tilde\xi$, after interpreting it as Fourier-Mukai
transform on $Y$, maps SLAGs on $X$ to holomorphic cycles on $X\dual$.
So the picture is roughly the following
$$\xymatrix{
X=(M,I)\ar @{>}[d]^{{\rm HK-rotation}}&\\
Y=(M,J)\ar @{>}[r]^{{\rm FM}=\xi}&Y=(M,J)\ar@{>}[d]^{{\rm HK-rotation}}\\
&X\dual=(M,-K)\\}
$$


\subsection{Large complex structure limit}

Mirror symmetry for Calabi-Yau manifolds suggests that the moduli space
of complex structures on one Calabi-Yau manifold should be canonically
isomorphic to the moduli space of K{\"a}hler structures on its dual. In
fact, this should literally only be true near  certain limit points.
The limit point for the complex structure is called large complex
structure limit and should correspond via mirror symmetry to the large
K{\"a}hler limit. Moreover, according to the SYZ version of
mirror symmetry, near the large complex structure limit the Calabi-Yau
manifold is a Lagrangian fibration and while approaching the limit the
Lagrangian fibres shrink to zero.

Let us consider a K3 surface $X$ together with a sequence of K{\"a}hler
classes $\omega_t$ such that the volume of $X$ with respect to
$\omega_t$ goes to infinity for $t\to\infty$. We will discuss in
particular cases, what happens to the mirror $X_t\dual$ of
$(X,\omega_t)$.
In order to incorporate the SYZ picture we will later concentrate on mirror
symmetry for elliptic K3 surfaces.

\medskip

There is a technical definition of the large complex structure limit
due to Morrison \cite{Morrisoncompl}. For a one-parameter family of K3
surfaces it goes as follows:

\begin{definition} A family of K3 surfaces $$\kx\dual\to D^*$$ over
the punctured disc $D^*:=\{z~|~0<|z|<1\}$ is a \emph{large complex structure
limit} if the monodromy operator $T$ on $H^2(X\dual,\IZ)$, where
$X\dual:=\kx\dual_t$
for some $t\ne0$, is maximally unipotent, i.e.\
$(T-1)^2\ne0$ but $(T-1)^3=0$, and $N:=\log(T)$ induces a filtration
$$W_0:={\rm Im}(N^2),~ W_1:={\rm Im}(N|_{\ker(N^2)}), ~W_2:={\rm
  Im}(N),~ {\rm and} ~W_3:=\ker(N^2),$$ such that
$\dim(W_0)=\dim(W_1)=1$ and $\dim(W_2)=2$.
\end{definition}

Geometrically such families arise as type III degenerations of K3
surfaces, as studied by Kulikov (cf.\ \cite{Friedman}).

\begin{remark} The weight filtration $W_*$ always satisfies $W_0\perp W_3$ (To see
this one use
$\langle (T-1)(a),b\rangle=-\langle T(a),(T-1)(b)\rangle$). In
particular, $W_0=W_1$ is spanned by an isotropic vector $v\in \Gamma$.

Conversely, any isotropic vector $v\in\Gamma$ together with an
additional
B-field $B_0\in\Gamma\setminus\IZ v$ defines a weight filtration
$W_0=W_1=\langle v\rangle$, $W_2=\langle v,B_0\rangle$, $W_3=v^\perp$
as above.
\end{remark}

Let us check that the mirror of a large K{\"a}hler limit 
is a large complex structure limit in the above sense.
We fix a decomposition $\Gamma=\Gamma'\oplus U'$ as before.
Now consider a K3 surface $X$ with a family of K{\"a}hler structures
$\omega_t\in\Gamma'_\IR$.

Recall that by Proposition \ref{msmap} the period of the mirror
$X_t\dual$ of $(X,\omega_t)$ is given by $$[\sigma_t\dual]
=[{\rm pr}(i\omega_t)+(\omega^2_t/2)v+v^*]\in\IP(\Gamma_\IC).$$
(For simplicity we assume $B=0$.) For $\omega_t^2\to\infty$ the
period point of $X\dual_t$  converges to $[v]$. The loop around the
large K{\"a}hler limit is given by $sB_0+i\omega_t$ with $s\in[0,1]$.
For $B_0\in\Gamma'$ the periods of the mirror for $s=0$ and $s=1$ are
thus given by $$i\omega_t+(\omega_t^2/2)v+v^*~{\rm resp.} ~i\omega_t-\langle\omega_t,B_0\rangle v+
(\omega_t^2/2)v+v^*+B_0-(B_0^2/2)v.$$

Thus, the induced monodromy is given by $T\in\OO(\Gamma)$
with $$T(v)=v, ~T(v^*)=v^*+B_0-(B_0^2/2)v,~{\rm and}~T(x)=x-\langle
B_0,x\rangle v~{\rm for}~ x\in\Gamma'.$$ Note that this is an  `internal
B-shift' by $B_0$ with respect to the decomposition
$\Gamma=\Gamma'\oplus U'$ of the type we have encountered already in
the proof of Proposition \ref{denseorb}. Hence, $T-1$ maps
$v\mapsto0$, $v^*\mapsto B_0-(B_0^2/2)v$, $x\mapsto-\langle
x,B_0\rangle v$ for $x\in\Gamma'$ and, therefore, $W_0=W_1$ is spanned by
$v$ and $W_2$ is spanned by $v$ and $B_0$.

This yields

\begin{proposition} Let us consider the mirror map induced
  by the decomposition $\Gamma=\Gamma'\oplus U'$.
  For any choice of $0\ne B_0\in\Gamma'$ the mirror of the
  large K{\"a}hler limit   $(X,\omega_t\in\Gamma'_\IR,B=0)$
  is a large complex structure limit. The choice of
  $B_0$ corresponds to choosing a component of the boundary divisor
  around which the monodromy is considered.\qed
\end{proposition}
In \cite{Dolgachev} the relation between the choice
of the decomposition $\Gamma=\Gamma'\oplus U'$ and the large complex
structure limit  is expressed by: `the choice of the isotropic
vector is the analog of MS1' (MS1: the choice of a boundary point with
maximally unipotent monodromy).

The SYZ conjecture can be incorporated into this picture without too much
trouble:
Let $\kx\dual\to D^*$ be the large complex structure limit obtained as
above endowed with the mirror K{\"a}hler structures $\omega\dual_t$. The
space $W_0$ is spanned by $v$.

\begin{proposition}
  The mirror $(\kx\dual_t,\omega\dual_t)$ of
  $(X,\omega_t)\in\Gamma'_\IR$ admits a
  SLAG fibration $\kx\dual_t\to\IP^1$ with fibre class $v$. The volume
  of the fibre converges to zero for $t\to0$.
\end{proposition}

\prf Recall that the holomorphic volume form
$\sigma$ is always chosen such that $\langle\im(\sigma),v\rangle=0$.
See the discussion at the beginning of the proof of Proposition
\ref{msmap}. Moreover, if an additional K{\"a}hler form $\omega$
is chosen one requires $\sigma\wedge\bar\sigma=2\omega^2$.

The class $v$ is isotropic, as was remarked before. It thus
satisfies a necessary condition for a fibre class.
Furthermore, using the formulae for the mirror one finds
that $\langle\omega\dual_t,v\rangle=\langle\im(\sigma\dual_t),v\rangle=0$,
where we use $\omega_t\in\Gamma'_\IR$. Thus, on the level of
cohomology the class $v$ is a SLAG.

Also note that
$\langle\re(\sigma_t\dual),v\rangle=\langle\re(\sigma_t),v\rangle^{-1}$.
Although the complex structure of $X$ does not change,
we have to rescale $\sigma$ in order to ensure
that $\sigma\bar\sigma=2\omega_t^2$.
The conditions on the choice of the real and imaginary part of
$\sigma_t$ imply that $\langle\re(\sigma_t),v\rangle\to\infty$. Hence,
$\langle\re(\sigma_t\dual),v\rangle\to 0$.
Therefore, the volume of the SLAG fibre class $v$ on $\kx\dual_t$
approaches
zero.

It thus remains to realize $v$ geometrically. Here one uses
hyperk{\"a}hler rotation as before. Indeed, rotating with respect to the
K{\"a}hler form $\omega\dual_t$ one finds a complex structure with
respect to which $v$ is of type $(1,1)$
and can thus be realized as the fibre class
of an elliptic fibration (modulo the action of the Weyl group).
\qed

\bigskip

It is tempting to apply this discussion to the case of elliptic K3
surfaces or to the case where mirror symmetry is described by
hyperk{\"a}hler rotation.
However, it seems impossible to follow the mirror to the large complex
structure limit and at the same time to be able to obtain the mirror
by a hyperk{\"a}hler rotation. This also supports the point of view
expressed in \ref{doubtsHKMS}, that mirror symmetry and hyperk{\"a}hler rotation
are only related to each other in a very restricted sense.

Let us recall the setting. We fix an elliptic K3 surface
$Y=(M,J)\to \IP^1$ with a section $\sigma_0$ and assume that the K{\"a}hler form
$\omega_J$ is of the form $(\alpha/2+1)f+\sigma_0$. The holomorphic
two-form $\sigma_J$ on $Y$ is chosen such that $\sigma_J\bar\sigma_J=2\omega_J^2$.
In this case the mirror of $X=(M,I)$ endowed with
$\omega_I=\im(\sigma_J)$
is $(M,-K)$.
In principal there are two ways one could try to combine the
passage of the mirror to the complex structure limit and the
interpretation of the mirror by hyperk{\"a}hler rotation.
First, one could fix the complex structure $I$, i.e.\ the K3 surface
$X=(M,I)$, and change the K{\"a}hler form $\omega_I$.

\medskip
\psset{unit=0.4cm}
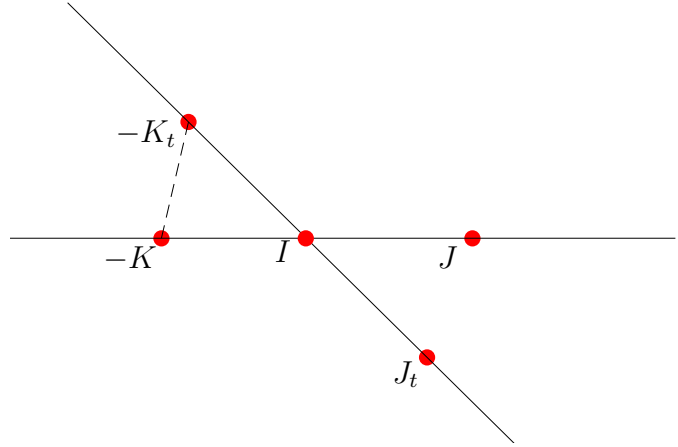
\hspace{3cm}\begin{pspicture*}(1,1)(23.100000,15.633333)
\psset{dotstyle=x, dotscale=2.0, linewidth=0.02}
\psdots[dotscale=2.5,dotstyle=*,linecolor=red](6.016667,7.800000)
\put(4.1,6.8){$-K$}
\psdots[dotscale=2.5,dotstyle=*,linecolor=red](16.350000,7.800000)
\put(15.2,6.8){$J$}
\psdots[dotscale=2.5,dotstyle=*,linecolor=red](6.916667,11.666667)
\put(4.5,11){$-K_t$}
\psdots[dotscale=2.5,dotstyle=*,linecolor=red](14.850000,3.833334)
\put(13.7,3){$J_t$}
\psline[linecolor=black](0.000000,7.799999)(23.100000,7.800001)
\psline[linecolor=black](2.899362,15.633333)(22.513263,-3.733333)
\psdots[dotscale=2.5,dotstyle=*,linecolor=red](10.816696,7.800000)
\put(9.8,7){$I$}
\psline[linestyle=dashed,linecolor=black](6.016667,7.800000)(6.916667,11.666667)
\end{pspicture*}
\medskip

However, in this case we keep $\sigma_I$. Rescaling is not permitted
if the K{\"a}hler form on $(M,J_t)$ is supposed to be of the form
$(\alpha/2+1)f+\sigma_0$. Thus,
$\langle\re(\sigma_I),f\rangle=\langle\omega_{J_t},f\rangle\equiv1$.
Hence, the volume of the SLAG fibres in the mirror does not converge
to zero. This yields a contradiction.

The second possibility is to freeze $Y=(M,J)$ and to change
$\omega_J=(\alpha/2+1)f+\sigma_0$ by considering
$(t^2\alpha/2+1)f+\sigma_0$ with $t\to \infty$.
In this case, we compute the mirror $X_t\dual$ of $X_t=(M,I_t)$
endowed with $t\omega_I$, which stays of type $(1,1)$ with respect
to the changing complex structure $I_t$.

\medskip

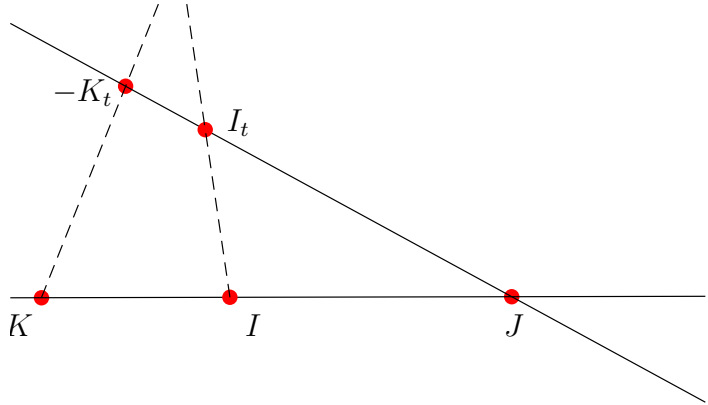
\psset{unit=0.8cm}
\hspace{3cm}\begin{pspicture*}(1,0)(12.550000,7.816667)
\psset{dotstyle=x, dotscale=2.0, linewidth=0.02}
\psdots[dotscale=2,dotstyle=*,linecolor=red](1.516667,2.933333)
\put(0.5,2.3){$-K$}
\psdots[dotscale=2,dotstyle=*,linecolor=red](9.333333,2.950000)
\put(9.2,2.3){$J$}
\psdots[dotscale=2,dotstyle=*,linecolor=red](2.916667,6.450000)
\put(1.7,6.2){$-K_t$}
\psline[linecolor=black](1.000000,2.932232)(12.550000,2.956859)
\psline[linecolor=black](1.000000,7.495455)(12.550000,1.195455)
\psline[linestyle=dashed,linecolor=black](1.516667,2.933333)(3.460743,7.816667)
\psdots[dotscale=2,dotstyle=*,linecolor=red](4.649949,2.940014)
\put(4.9,2.3){$I$}
\psdots[dotscale=2,dotstyle=*,linecolor=red](4.239707,5.728342)
\put(4.6,5.7){$I_t$}
\psline[linestyle=dashed,linecolor=black](4.649949,2.940014)(3.932454,7.816667)
\end{pspicture*}

\medskip

In this case, one finds that the complex structures $-K_t$ and $I_t$
converge. But as before, the fibre volume
$\langle\re(\sigma_{I_t},f\rangle
=\langle (t^2\alpha/2+1)f+\sigma_0,f\rangle$ does not converge to zero.

\bigskip

The large complex structure limit was studied by Gross and Wilson on a
much deeper level in \cite{GW}. They studied the corresponding
sequence of K{\"a}hler metrics (or at least a very good approximation
of those) and could indeed show that the fibres of the elliptic
fibration are shrunk to zero. Moreover, they managed to describe the
limit metric on the base sphere $S^2=\IP^1$.


\bigskip

\noindent

\end{document}